\documentstyle[pb-diagram]{article}

\newcommand{\es}{\emptyset}

\newcommand{\ba}{\begin{array}}
\newcommand{\ea}{\end{array}}

\newcommand{\en}{en}

\newtheorem{theorem}{THEOREM}[section]
\newtheorem{lemma}{LEMMA}[section]
\newtheorem{corollary}{COROLLARY}[section]
\newtheorem{definition}{DEFINITION}[section]
\newtheorem{proposition}{PROPOSITION}[section]

\newcommand{\be}{\begin{enumerate}}
\newcommand{\ee}{\end{enumerate}}
\newcommand{\bi}{\begin{itemize}}
\newcommand{\ei}{\end{itemize}}
\newcommand{\bd}{\begin{description}}
\newcommand{\ed}{\end{description}}

\newcommand{\et}{\wedge}

\newcommand{\vel}{\vee}

\newcommand{\imp}{\rightarrow}

\newcommand{\beq}{\begin{eqnarray*}}
\newcommand{\eeq}{\end{eqnarray*}}
\newcommand{\seq}{\Rightarrow}

\author{ {F.Parlamento, F.Previale }
\\Department of Mathematics,  Computer Science and Physics
\\University of Udine,  via  Delle Scienze 206, 33100 Udine, Italy.
\\Department of Mathematics
\\University of Turin, via Carlo Alberto 10, 10123 Torino, Italy
\\e-mail: {\em franco.parlamento$@$uniud.it},  {\em flavio.previale$@$unito.it}
}
 \title{The Cut Elimination and Nonlengthening   Property for the  Sequent Calculus with Equality\thanks{Work  supported by funds PRIN-MIUR of Italy, Grant "Logica, Modelli, Insiemi"  and presented to the Logic Colloquium 2016(Leeds)} }
\date{20/03/2017}

\begin{document}
\maketitle



\begin{abstract}
We show how Leibnitz's indiscernibility principle and Gentzen's
 original work lead to extensions of the sequent calculus to  first order logic with equality and investigate the cut elimination property. Furthermore we discuss and improve  the nonlengthening property
of Lifschitz and Orevkov in \cite{L68} and \cite{O69}.
\end{abstract}

\section{Introduction}

 The most  common way  of treating equality in sequent calculus is to add to Gentzen's system
appropriate  sequents  with which   derivations can start, beside the logical axioms of the form $F\seq F$ (see for example \cite{C77}, \cite{T87}, \cite{TS96})  and \cite{Ga86}).
In this  way equality is considered and treated as a mathematical relation  subject to specific axioms.
For such kind of calculi Gentzen's cut elimination theorem can hold at most in a weakened form: every derivation can be transformed into one 
which contains only cuts whose cut formula is an equality. That doesn't allow to obtain directly the wealth of applications that
 full cut elimination  has, such as the conservativity  of first order logic  with equality over first order logic  without equality, 
or the disjunction and existence property for intuitionistic logic with equality.
As shown in \cite{NvP98},  the initial sequents  that concern equality can be replaced by nonlogical rules in order to obtain sequent calculi for which all the structural rules, including the cut rule, 
are admissible.
 However the above  nonlogical  rules   can  eliminate equalities, so that obtaining the mentioned applications of cut elimination is not entirely straightforward
 (see \cite{NvP98}, \cite{NvP01} or \cite{TS00},   for the additional work required to obtain the conservativeness of
first order logic  with equality over first order logic without equality).
Our purpose is to overcome this  difficulty by introducing  a sequent calculus for which full cut elimination holds and none of the  rules, other than the cut rule,  eliminates equalities or  any other logical constant, and
moreover, to remain as close as possible to Gentzen's system, retains the  separation between structural and logical rules.

 For notational semplicity and to add  evidence to  the logical definability of equality, throughout this introduction and the next section, we will restrict attention to intuitionistic logic.

To begin with,  we observe that   equality can be regarded as a logical   constant, defined,  according to Leibnitz's indiscernibility principle, by letting $a=b$ to mean $\forall X(X(a)\leftrightarrow X(b))$.
 However, in the framework of intuitionistic second order logic,  thanks to the  rules for $\forall $ and $\imp$,
 $\forall X(X(a)\leftrightarrow X(b))$ is equivalent to  $\forall X(X(a)\rightarrow X(b))$.
In fact $\forall X(X(b)\rightarrow X(a))$ can be deduced from $\forall X(X(a)\rightarrow X(b))$ by instantiating the bound predicate variable $X$ by the lambda term $\lambda v (Z(v)\imp Z(a))$, where $Z$ is a free
predicate variable,  so as to obtain  $(Z(a)\imp Z(a))\imp (Z(b)\imp Z(a))$. Then, given the deducibilty, by $\imp$-introduction,  of $Z(a)\imp Z(a)$, an $\imp$-elimination followed by a $\forall$-introduction yields $\forall X(X(b)\rightarrow X(a))$ as desired.
Thus as the  definition of $a=b$ we can simply take  $\forall X (X(a)\imp X(b))$. On that respect equality 
stands on a par with the definition of $\et, \vel, \neg, \exists$ 
in terms of universal quantification   and implication, spelled out, for example, in   \cite{Pr65}) pg. 67.
Given this   definition of $=$, from the rules of Gentzen's  sequent calculus for $\forall$ and $\imp$,  the following
left and right introduction rules for $=$:

\[
\ba{clccc}
\Lambda \seq F\{v/r\}~~~~~\Gamma, F\{v/s\} \seq \Delta& \vbox to 0pt{\hbox{$=^{(2)}\seq $}}&~~~~&\Gamma, Z(r)\seq Z(s)&\vbox to 0pt{\hbox{$\seq =^{(2)} $}}\\
\cline{1-1} \cline{4-4}
\Lambda, \Gamma, r=s\seq \Delta&&&\Gamma \seq r=s
\ea
\]can be derived. Here and in the following  $F\{v/r\}$ ($F\{v/s\}$) denotes the result of the simultaneous replacement in $F$ of all the occurrences of the free object variable   $v$  by $r$ ($s$) and $Z$ is a free predicate variable that  does not occur in   $\Gamma$,  and $|\Delta|\leq 1$. 
Conversely the sequents $r=s\seq \forall X( X(r)\imp X(s))$ and $\forall X( X(r)\imp X(s))\seq r=s$ are derivable by using the rules 
 $=^{(2)}\seq $ and  $\seq =^{(2)} $
(the details of such derivations as well as  of a few others to  be  mentioned in this  Introduction are provided in the next section  of the paper).
Granted Leibniz's definition of equality, we can therefore claim that the second order version of $LJ$ supplemented by the rules
$=^{(2)}\seq $ and  $\seq =^{(2)} $, that we denote by $LJ^{(2)=}$, is an adequate  sequent calculus to deal with equality in second order logic. 
 The  right introduction rule $\seq =^{(2)}$ turns out to be equivalent to
the {\em Reflexivity Axiom} $\seq r=r$.
Thus a sequent calculus for first order logic with equality  can  be obtained   from $LJ^{(2)=}$ by  replacing  $\seq  =^{(2)}$ by the Reflexivity Axiom  and        requiring  that all the formulae and terms involved  be  first order formulae and terms.

We will denote by $=\seq$ and $\seq =$ the rule and axiom obtained in that way, and by $LJ^{(1)=}$
 the sequent calculus that is obtained by adding them to Gentzen's  $LJ$.
 $LJ^{(1)=}$ shares with $LJ$ the distinction between structural and logical rules and the latter introduce a logical constant as the outermost symbol of exactly one formula,  the so called   {\em principal formula}, while the other formulae in the conclusion are those  present in a determined  position in the premiss or premisses and are independent from the principal one.
 
  As we will see, full cut elimination holds for $LJ^{(1)=}$, 
however  $LJ^{(1)=}$ is far from being a satisfactory sequent calculus for first order logic with equality,
since the application of the rule $=\seq$ eliminates  the formulae $F\{v/r\}$ and $F\{v/s\}$,
hence all the logical constants they may contain.

Our task is therefore to find a calculus equivalent to $LJ^{(1)=}$ 
in which the cut rule is eliminable and all the other rules do not eliminate occurrences of logical constants.
Following the lines of Gentzen's transition from the axiomatic systems to natural deduction and then to the sequent calculus  (see  \cite{vP12} and \cite{vP14} for a  detailed historical reconstruction),  
our  starting point will be the following natural deduction elimination rule for $=$:

\[
\ba{cl}
F\{v/r\}~~~~~r=s&\vbox to 0pt{\hbox{$=_1^N $}}\\
\cline{1-1}
F\{v/s\}&
\ea
\]
together with the rule for its introduction, namely  the zero premisses  reflexivity rule $\overline{r=r}$, that correspond to the substitutivity axiom  $\forall x\forall y (x=y \imp (F\{v/x\} \imp F\{v/y\})$ and to the  reflexivity axiom $\forall x (x=x)$.

Let $NJ^=$ be the natural deduction system obtained by adding to $NJ$ the above introduction and elimination rules for $=$.
We will pick  the right and left introduction rules for $=$ to be added to $LJ$ so  as to obtain a Gentzen-style sequent calculus equivalent to $NJ^=$, namely such  that
a formula $G$ is deducible from (assumptions that are listed in) $\Sigma$ if and only if $\Sigma \seq G$ is derivable in the calculus.
Since $r=r$ is deducible from the empty $\Sigma$, the most obvious corresponding choice is to add to $LJ$, 
as the (zero premisses) right introduction rule, 
the {\em Reflexivity Axiom} $\seq r=r$ 
(already denoted by  $\seq =$). Considering  the correspondence between the natural deduction elimination rules and the left introduction rules of the sequent calculus, particularly those concerning the existential quantifier,
it is quite natural to make correspond to $=_1^N$ the following  rule:

\[
\ba{cl}
\Gamma\seq F\{v/r\}&\vbox to 0pt{\hbox{$=_1 $}}\\
\cline{1-1}
\Gamma, r=s\seq F\{v/s\}&
\ea
\]
Actually, the most direct transformation in sequent terms of $=_1^N$ is the following rule, clearly equivalent to $=_1$ over the structural rules:

\[
\ba{cl}
\Gamma\seq F\{v/r\}~~~~~\Lambda \seq r=s&\vbox to 0pt{\hbox{$~CNG$}}\\
\cline{1-1}
\Gamma, \Lambda\seq F\{v/s\}
\ea
\]
($CNG$ for {\em congruence}), that will play a crucial role in the sequel. 
Given the equivalence between $=_1$ and $CNG$  it is straightforward that, if we let $LJ_\circ^=$ be the calculus obtained by adding $=_1$ and $=_2$ to $LJ$, then $LJ_\circ^=$ is equivalent to $NJ^=$.
However, cut elimination  for $LJ_\circ^=$  does not hold.
For example the following derivable sequent $a=c,b=c \seq a=b$  cannot have any cut free derivation, if we adopt only $\seq =$ and $=_1$.

In order to have cut elimination we have to add also the following rule, obtained by replacing in $=_1$, $r=s$ by its symmetric $s=r$:

\[
\ba{cl}
\Gamma\seq F\{v/r\} &\vbox to 0pt{\hbox{$=_2 $}}\\
\cline{1-1}
\Gamma, s=r\seq F\{v/s\}
\ea
\]

corresponding to the following other natural deduction elimination rule for $=$:

\[
\ba{cl}
F\{v/r\}~~~~~s=r&\vbox to 0pt{\hbox{$=_2^N $}}\\
\cline{1-1}
F\{v/s\}&
\ea
\]

Letting $LJ^=$ be the result of adding to $LJ$ both $=_1$ and $=_2$
we will provide a very simple proof that cut elimination holds for $LJ^=$, based on the admissibility in the cut-free part of $LJ^=$ of the rule  $CNG$ introduced above (which by itself would seem of scarce interest for the sequent calculus, since its application eliminates equalities). $LJ^=$ and $LJ^{(1)=}$ are equivalent and  $=_1 $ and $=_2 $ are   derivable in $LJ^{(1)=}$ without using the cut rule. Hence 
cut elimination for $LJ^{(1)=}$ follows as an immediate consequence of cut elimination for $LJ^=$.

In the light of the derivability   of $=_1$ and $=_2$ in $LJ^{(1)=}$, it is quite natural to consider also
the following rules:
$=_1^l$ and   $=_2^l$:

\[
\ba{clccl}
\Gamma, F\{v/r\}\seq \Delta&\vbox to 0pt{\hbox{$=_1^l$ }}&~~~\mbox{and}~~~&\Gamma, F\{v/r\}\seq \Delta&\vbox to 0pt{\hbox{$=_2^l$ }}\\
\cline{1-1}\cline{4-4}
\Gamma, F\{v/s\}, r=s\seq \Delta&&&\Gamma, F\{v/s\}, s=r\seq \Delta
\ea
\]
The  four equality rules $=_1$, $=_2$, $=_1^l$ and $=_2^l$ turn out to be equivalent to $=\seq$, hence to each other, over the structural rules and $\seq =$.

We will show that cut elimination holds for the systems,  to be denoted by $LJ_1^=$ and $LJ_2^=$, that are obtained from  $LJ$  by adding $\seq =$,  $=_1$ and $=_1^l$ or    $\seq =$, $=_2$ and $=_2^l$.
In fact we will show that the rule $=_2$ is admissible  in $LJ^=_1$ deprived of the cut rule and, similarly,  that $=_1$ is  admissible in  $LJ^=_2$ 
deprived of the cut rule,  so that cut elimination for both systems
follows form cut elimination for $LJ^=$.

Despite the similarity of   the
pair of rules $=_1$ and $=_2$ and the pair $=_1^l$ and $=_2^l$ with respect to $LJ^{(1)=}$,
  the system obtained from $LJ$ by adding   $\seq =$ and  both $=_1^l$ and  $=_2^l$ does not satisfy cut elimination.
That turns out to be the case  also for the systems that are obtained from  $LJ$  by adding  $\seq =$ together with $=_1$ and $=_2^l$ or $\seq =$ together with $=_1^l$ and $=_2$.

Furthermore we will show   that if all the four equality rules are adopted, then we obtain a system $LJ^=_{12}$ for which
cut elimination holds also if  their application is  required to be $\prec$-nonlengthening, with respect to any binary antisymmetric relation on terms $\prec$. We recall  from  \cite{L68}, that an equality-inference as represented above is said to be $\prec$- nonlengthening if $s\not\prec r$. Actually we  will show that cut elimination holds for the system in which all the equality-inferences are required to be $\prec$-nonlengthening and all the $=_1$ and $=_1^l$-inferences are required to be $\prec$-shorthening, namely to satisfy  the stronger condition $r\prec s$. Alternatively we can require that all the equality-inferences be $\prec$-nonlengthening and all the $=_2$ and $=_2^l$-inferences be $\prec$-shorthening.

All the above results hold  without any essential change for  the classical version of the calculi considered, in particular for the classical version $LK^=$ of $LJ^=$.

The  union $LK^=_{12}$ of the systems $LK^=_1$ and $LK^=_2$ is equivalent, on the  ground of the exchange  and contraction rules only, to the system $G^e$
 in \cite{L68}, that was motivated by the calculus free of structural rules introduced in \cite{K63}, for efficient proof search  purposes.
Therefore we have a proof that, as announced in \cite{L68}, $G^e$ satisfies cut elimination.
Finally, improving the result stated in \cite{O69},   for any antisymmetric relation $\prec$ on terms, we will show that 
any derivation in $LK^=_{12}$,  can be transformed into a cut-free derivation in the  same system of its endsequent, whose equality inferences are $\prec$-nonlengthening or $\prec$- shortening, as explained above for the intuitionistic case.

\subsection{Basic Derivations}

Having defined   $r=s$ as $\forall X (X(r)\imp X(s)$, the conclusion of the rule $=^{(2)}\seq $, namely

\[
\ba{c}
\Lambda \seq F\{v/r\}~~~~~~~~\Gamma, F\{v/s\} \seq \Delta\\
\cline{1-1}
\Lambda, \Gamma, r=s \seq \Delta
\ea
\]
can be derived from its premisses by 
applying first the left introduction rule for $\imp$ and then the second order left introduction rule for $\forall$, while the conclusion of $\seq =^{(2)}$, namely

\[
\ba{c}
\Gamma, Z(r)\seq Z(s)\\
\cline{1-1}
\Gamma\seq r=s
\ea
\] can be derived from its premiss by applying first  the right introduction rule for $\imp$ and  then the second order right introduction rule for $\forall$.

Conversely the sequents $r=s\seq \forall X(X(r)\imp X(s)$ and $\forall X(X(s)\imp X(r)\seq r=s$ can be derived by means of $=^{(2)}\seq$ and $\seq=^{(2)} $ as follows:
 \[
 \ba{cccl}
 Z(r)\seq Z(r)~~~~~Z(s)\seq Z(s)&\vbox to 0pt{\hbox{$=^{(2)}\seq $}}~~~~~~&Z(r)\seq Z(r)~~~~~Z(s)\seq Z(s)&\\
 \cline{1-1}\cline{3-3}
 Z(r), r=s\seq Z(s)&& Z(r)\imp Z(s), Z(r)\seq Z(s)&\\
 \cline{1-1}\cline{3-3}
 r=s\seq Z(r)\imp Z(s)&&\forall X(X(r)\imp X(s)), Z(r)\seq Z(s)&\vbox to 0pt{\hbox{$\seq  =^{(2)}$}}\\
 \cline{1-1}\cline{3-3}
 r=s\seq \forall X (X(r)\imp X(s))&&\forall X(X(r)\imp X(s))\seq r=s&
 \ea
 \]
where we have omitted the applications of the exchange rule, as we will do throughout the paper.

$\seq r=r$ is immediately derived by  $\seq =^{(2)}$ applied to the logical axiom $Z(r)\seq Z(r)$ and, 
 conversely, $\seq =^{(2)}$, can be derived from $\seq r=r$, by using the cut rule,   as follows:

\[
\ba{ccc}
\Gamma, Z(r)\seq Z(s)&~~~~&\seq r=r~~~~~r=s\seq r=s\\
\cline{1-1}\cline{3-3}
\Gamma\seq Z(r)\imp Z(s)&&r=r\imp r=s\seq r=s\\
\cline{1-1}\cline{3-3}
\Gamma\seq \forall X (X(r)\imp X(s))&&\forall X (X(r)\imp X(s))\seq r=s\\
\cline{1-3}
\multicolumn{3}{c}{\Gamma \seq r=s}
\ea
\]

\

Concerning the equivalence between $LJ^=$ and $LJ^{(1)=}$ we note that both $=_1$ and $=_2$ are derivable from $=\seq$, without using the cut rule.
In fact  we have the following derivations of $=_1$ and $=_2$ respectively:

\[
\ba{cccl}
\Gamma\seq F\{v/r\}&&F\{v/s\}\seq F\{v/s\}&\vbox to 0pt{\hbox{$=\seq $}}\\
\cline{1-3}
\multicolumn{3}{c}{\Gamma, r=s \seq F\{v/s\}}&
\ea
\]

\[
\ba{cccccl}

&&\Gamma\seq F\{v/r\}&&F\{v/s\}\seq F\{v/s\}&\vbox to 0pt{\hbox{$=\seq $}}\\
 \cline{3-5}
\seq s=s&&\multicolumn{3}{c}{\Gamma, r=s\seq F\{v/s\}}&\vbox to 0pt{\hbox{$=\seq $}}\\
\cline{1-5}
\multicolumn{5}{c}{\Gamma, s=r \seq F\{v/s\}}
\ea
\]
where the last inference is a correct application of $=\seq$ in  which the place of $F$ is taken by  $v=s$.
Conversely,  by using the cut rule,  $=\seq$   can be derived from $=_1$ and also from $=_2$ as follows:

\[
\ba{clcc}
\Lambda  \seq F\{v/r\}&\vbox to 0pt{\hbox{$=_1$}}&~~~~~~~&\\
\cline{1-1}
\Lambda, r=s \seq F\{v/s\}&&~~~~~\Gamma, F\{v/s\}\seq \Delta\\
\cline{1-4}
\multicolumn{4}{c}{\Lambda, \Gamma,r=s\seq \Delta}
\ea
\]

\[
\ba{clclc}
\seq s=s&\vbox to 0pt{\hbox{$=_2$}}~~~&\Lambda \seq F\{v/r\}&\vbox to 0pt{\hbox{$=_2$}}~~~&\\
\cline{1-1}\cline{3-3}
r=s\seq s=r &&\Lambda, s=r\seq F\{v/s\}&&\\
\cline{1-3}
\multicolumn{3}{c}{~~~~\Lambda, r=s \seq F\{v/s\}}& \Gamma, F\{v/s\} \seq \Delta\\
\cline{2-4}
&\multicolumn{3}{c}{\Lambda, \Gamma, r=s \seq \Delta}

\ea
\]

Therefore it would suffice to add $\seq =$ and $=_1$ or   $\seq =$ and $=_2$ to $LJ$ in order to have a system equivalent to $LJ^{(1)=}$.

As for  $=_1^l$ and $=_2^l, $ namely:

\[
\ba{clccl}
\Gamma, F\{v/r\}\seq \Delta&&~~~\mbox{and}~~~&\Gamma, F\{v/r\}\seq \Delta&\\
\cline{1-1}\cline{4-4}
\Gamma, F\{v/s\}, r=s\seq \Delta&&&\Gamma, F\{v/s\}, s=r\seq \Delta
\ea
\] we have the following derivations  in $LJ^{(1)=}$:

\[
\ba{cccl}
&&F\{v/s\}\seq F\{v/s\}~~~~\Gamma, F\{v/r\} \seq \Delta&\vbox to 0pt{\hbox{$=\seq $}}\\
\cline{3-3}
\seq r=r&~~~&\Gamma, F\{v/s\}, s=r \seq \Delta&\vbox to 0pt{\hbox{$=\seq $}}\\
\cline{1-3}
\multicolumn{3}{c}{\Gamma, F\{v/s\}, r=s \seq \Delta}
\ea
\]and

\[
\ba{cl}
F\{v/s\}\seq F\{v/s\}~~~~\Gamma, F\{v/r\} \seq \Delta& \vbox to 0pt{\hbox{$=\seq$ }}\\
\cline{1-1}
\Gamma, F\{v/s\}, s=r \seq \Delta
\ea
\]

Conversely the  rule  $=\seq$ can be derived from  $=_2^l$ or $\seq =$ and $=_1^l$  as follows:

\[
\ba{cccl}
&&\Lambda, F\{v/s\}\seq \Delta&\vbox to 0pt{\hbox{$=_2^l$}}\\
\cline{3-3}
\Gamma\seq F\{v/r\}&~~~&\Lambda, F\{v/r\}, r=s\seq \Delta&\\
\cline{1-3}
\multicolumn{3}{c}{\Gamma, \Lambda, r=s  \seq \Delta}
\ea
\]

\[
\ba{cccl}
&&\Gamma, F\{v/s\}\seq \Delta&\vbox to 0pt{\hbox{$=_1^l$}}\\
\cline{3-3} 
&&\Lambda, F\{v/r\}, s=r \seq \Delta&\vbox to 0pt{\hbox{$=_1^l$}}\\
\cline{3-3}
&\seq s=s&\Lambda, F\{v/r\}, s=s, r=s \seq \Delta&\\
\cline{2-3}
\Gamma\seq F\{v/r\}&&\Lambda, F\{v/r\}, r=s\seq \Delta&\\
\cline{1-3}
\multicolumn{3}{c}{\Gamma, \Lambda, r=s  \seq \Delta}

\ea
\]

Therefore $=_1$, $=_2$, $=_1^l$ and $=_2^l$ are all equivalent to $=\seq$, and therefore to each other,  over $\seq =$ and the structural rules.


 \subsection{Transformation  of derivations into separated form}
 In the  following  $LJ$ and $LK$  will denote  the sequent calculi introduced by Gentzen in \cite{G35},  except that,  as in \cite{T87},  in the left introduction 
 rule $\forall\seq$ for $\forall$ and in the right introduction rule 
  $\seq \exists$ for $\exists$ the {\em free object variable} is replaced by an arbitrary term.
  Clearly what  has been said so far about $LJ$ holds for $LK$ as well, with the obvious changes in the presentation of the rules, needed to allow the possible presence of more than one formula in the succedent of the sequents.
 
 
$LJ^=$ and $LK^{=}$ are obtained   by adding to $LJ$ and $LK$, the {\em equality rules} $=_1$ and $=_2$, namely

\[
\ba{clccl}
\Gamma\seq \Delta, F\{v/r\}&\vbox to 0pt{\hbox{ $=_1$}}~~~\mbox{and}~~~~&\Gamma\seq \Delta, F\{v/r\}&\vbox to 0pt{\hbox{ $=_2$}}\\
\cline{1-1}\cline{3-3}
\Gamma, r=s\seq \Delta, F\{v/s\}&&\Gamma, s=r\seq \Delta, F\{v/s\}
\ea
\]
where    $v$ is a free object variable that occurs neither in $r$ nor in $s$ and, in the case of $LJ^=$, $\Delta=\es$. Notice  that  the requirement on $v$ is not restrictive  since $F\{v/r\}$ and $F\{v/s\}$  can always be represented as
$(F\{v/v'\})\{v'/r\}$ and  $(F\{v/v'\})\{v'/s\}$ for any  $v'$ that is  new to $F$, $r$ and $s$. If $v$ does not occur in $F$,  $=_1$ and $=_2$ reduce to a   left weakening, 
introducing $r=s$,  and are said to be {\em trivial}. $r=s$  in the presentation of $=_1$  and $s=r$   in the presentation of  $=_2$,  will be called the {\em operating equality}, while  
$F$  will be called the {\em changing  formula}  (in the representation)  of $=_1$ and $=_2$ .

\

At the  purely equational level, $LJ^=$ and $LK^=$ are equivalent, namely  a sequent $\Gamma\seq F$ is derivable in $LJ^=$ without applications of  logical rules, if and only if it is derivable in $LK^=$, without applications of  logical rules. In fact a straightforward induction on the height of derivations establishes the following:

\begin{proposition}\label{int}

If  a sequent $\Gamma\seq \Delta$ is derivable in   $LK^{=}$ without applications of logical rules, then there is a formula $F$ in $\Delta$ such  that $\Gamma\seq F$ has a derivation  without applications of logical rules, that contains only sequents with exactly one formula in the succedent. In particular $\Gamma\seq$ is not derivable in $LK^=$ without applications of  logical rules.
\end{proposition}

Proposition \ref{int} motivates the following definition:

\begin{definition}
$EQ$ is the calculus acting on  sequents with one formula in the succedent,   having the logical axioms $F\seq F$,   the reflexivity axioms  
$\seq t=t$;  the {\em weak left structural rules of weakening, exchange and contraction}:

\[
\ba{ccccc}
\Gamma\seq H&~~~&\Gamma_1,F,G,\Gamma_2\seq H&~~~&\Gamma,F, F\seq H\\
\cline{1-1}\cline{3-3}\cline{5-5}
\Gamma, F\seq H&&\Gamma_1, G, F, \Gamma_2\seq H&~~~&\Gamma,F\seq H
\ea
\]the {\em  cut rule}:
\[
\ba{c}
\Gamma \seq F~~~~\Lambda, F \seq H\\
\cline{1-1}
\Gamma,\Lambda  \seq H
\ea
\]and the   equality left  introduction 
rules $=_1$ and $=_2$:

\[
\ba{ccc}
\Gamma\seq F\{v/r\}&~~~~~~~&\Gamma\seq F\{v/r\}\\
\cline{1-1}
\cline{3-3}
\Gamma, r=s \seq F\{v/s\}&~~~~~~~&\Gamma, s=r\seq F\{v/s\}\
\ea
\].

\end{definition}

Our proof of cut elimination for $LJ^=$ and $LK^=$ will  split into two parts. First we show that every derivation can be tranformed into one that consists of derivations in $EQ$ followed by applications of weak structural rules, namely structural rules different from the cut rule, and logical rules only, and then that cut elimination holds for  $EQ$.

\begin{definition}

A derivation in $LJ^=$ or $LK^{=}$ is said to be {\em separated} if it consists of derivations in $EQ$, followed by applications (possibly none)  of logical and weak structural rules (both  left and  right)
\end{definition}

In order to prove that every derivation can be transformed into a separated derivation of its endsequent, we prove first that, thanks to the cut rule, the equality rules can be derived from  their special case in which the formula that they transform is atomic, and then that  the cut  rule is admissible over its restriction to atomic cut formulae.



\begin{lemma}\label{basic-atomic}
The sequents of the following form:
\bi
\item[a)] $F\{v/r\}, r=s\seq F\{v/s\}$
\item[b)] $F\{v/r\}, s=r\seq F\{v/s\}$
\ei
have derivations whose equality inferences are atomic, manely have the form

\[
\ba{ccc}
\Gamma\seq \Delta, A\{v/r\}&~~~~~&\Gamma \seq \Delta, A\{v/r\}\\
\cline{1-1}\cline{3-3}
\Gamma, r=s \seq \Delta, A\{v/s\}&~~~~~&\Gamma, s=r \seq \Delta, A\{v/s\}
\ea
\]
where $A$ is required to be an atomic formula.
\end{lemma}

 {\bf Proof} We proceed by induction on the degree of $F$. If $F$ is atomic, for $a)$ it suffices to consider
\[
\ba{cl}
F\{v/r\}\seq F\{v/r\}&\vbox to 0pt{\hbox{$=_1 $}}\\
\cline{1-1}
F\{v/r\}, r=s \seq F\{v/s\}
\ea
\]
As for $b)$ it suffices to replace $=_1$ by $=_2$.

If $F$ is $\neg G$, to establish $a)$, we apply  the  induction hypothesis $b)$  to  $G$,  according to which 
there is a derivation ${\cal D}$,  whose  equality inferences are atomic, of $G\{v/s\}, r=s\seq G\{v/r\}$. Then the following is the desired derivation:

\[
\ba{c}
{\cal D}\\
G\{v/s\}, r=s\seq G\{v/r\}\\
\cline{1-1}
G\{v/s\}, r=s,  \neg G\{v/r\}\seq \\
\cline{1-1}
 \neg G\{v/r\}, r=s \seq \neg G\{v/s\}
\ea
\]
$b)$ is established in a similar way by using the induction hypothesis $a)$ applied to $G$.

If $F$ is $G\imp H$, to establish $a)$ we apply the induction hypothesis $b)$ to $G$ and $a)$ to $H$ according to which there
are derivations ${\cal D}$ and ${\cal E}$ of $G\{v/s\}, r=s\seq G\{v/r\}$ and $H\{v/r\}, r=s\seq H\{v/s\}$ respectively,
whose equality inferences are atomic.
Then the following derivation establishes $a)$ for $F$:

\[
\ba{ccc}
{\cal D}&~~~~~~~&{\cal E}\\
G\{v/s\}, r=s\seq G\{v/r\}&&H\{v/r\}, r=s\seq H\{v/s\}\\
\cline{1-3}
\multicolumn{3}{c}{G\{v/s\}, r=s, r=s, G\{v/r\}\imp H\{v/r\} \seq H\{v/s\}}\\
\cline{1-3}
\multicolumn{3}{c}{G\{v/s\}, r=s,  G\{v/r\}\imp H\{v/r\} \seq H\{v/s\}}\\
\cline{1-3}
\multicolumn{3}{c}{      G\{v/r\}\imp H\{v/r\}, r=s \seq G\{v/s\} \imp H\{v/s\}    }
\ea
\]
$b)$ for $F$ is established in a similar way,  except that we have to use the induction hypothesis $a)$ on $G$ and $b)$ on $H$.

If $F$ is $\forall x G$, we let $u$ be any parameter not occurring in $G, r, s$ and apply the induction hypothesis $a)$ to $G\{x/u\}$ to obtain a  derivation ${\cal D}$, whose equality inferences are atomic, of $G\{v/r,x/u\}, r=s \seq G\{v/s,x/u\}$. Then the following derivation  establishes $a)$ for $F$:

\[
\ba{c}
G\{v/r,x/u\}, r=s \seq G\{v/s, x/u\}\\
\cline{1-1}
\forall x G\{v/r\}, r=s \seq G\{v/s, x/u\}\\
\cline{1-1}
\forall x G\{v/r\}, r=s \seq \forall x G\{v/s\}
\ea
\]
$b)$ for $F$ is established in the same way except that we use the induction hypothesis $b)$, rather than $a)$, on $G\{x/u\}$.

The other cases are similar and we omit the details $\Box$

\begin{proposition}\label{atomic}

Any  non atomic equality inference in a given derivation in $LJ^=$ or $LK^=$ can be replaced  by a cut between its premiss
and the endsequent of a derivation that uses only atomic equality inferences.  In particular any derivable sequent in $LJ^=$ or $LK^=$ has a derivation whose equality inferences are all atomic.
\end{proposition}

{\bf Proof} A non atomic $=_1$-inference of the form:
\[
\ba{c}
\Gamma \seq \Delta, F\{v/r\}\\
\cline{1-1}
\Gamma, r=s \seq \Delta, F\{v/s\}
\ea
\]
can be replaced by:

\[
\ba{ccc}
&~~~~~~~~~~&{\cal D}\\
\Gamma\seq \Delta, F\{v/r\}&& F\{v/r\}, r=s \seq F\{v/s\}\\
\cline{1-3}
\multicolumn{3}{c}{     \Gamma, r=s \seq \Delta, F\{v/s\}}
\ea
\]where ${\cal D}$ is the  derivation containing only atomic equality-inferences of Lemma \ref{basic-atomic} $a)$ for $F$.
A non atomic $=_2$-inference is eliminated in a similar way using  Lemma \ref{basic-atomic} $b)$. $\Box$

\

{\bf Notation} In the following $A$  will always denote  an atomic formula and 
$\Gamma\sharp F$ will denote any sequence of formulae from which $\Gamma$ can be obtained by eliminating any number, possibly none, of occurrences of $F$.




\begin{proposition}\label{atomic-cut}
If $\Gamma \seq \Delta \sharp F$ and $\Lambda \sharp F \seq \Theta$ have  derivations in $LJ^=$ or $LK^=$ whose equality and cut-inferences are  atomic, then also $\Gamma,\Lambda \seq \Delta,\Theta $ has a derivation in the same system whose equality and cut-inferences are atomic.
\end{proposition}

{\bf Proof}  Let  ${\cal D}$ and ${\cal E}$ be   derivations of $\Gamma\seq \Delta\sharp F$ and $\Lambda \sharp F \seq \Delta$ whose equality and cut-inferences are atomic.
If $\Delta\sharp F$ coincides with $\Delta$ or $\Lambda\sharp F$ coincides  with $\Lambda$,  then the desired  derivation of  $\Gamma, \Lambda \seq \Delta, \Theta$ can be simply obtained by applying some weakenings to the end sequent $\Gamma \seq \Delta$ of ${\cal D}$ or some exchanges and weakenings to the endsequent $\Lambda \seq \Theta$ of ${\cal E}$. 
 We can therefore assume that in $\Delta \sharp F$ there are occurrences of $F$ that are not listed in $\Delta$ and similarly for 
$\Lambda  \sharp F$.
 If $F$ occurs in $\Delta$, then from $\Gamma\seq \Delta\sharp F$ we can derive $\Gamma\seq \Delta$ by means of exchanges and contractions, and from $\Gamma\seq \Delta$ we can then derive $\Gamma,\Lambda \seq \Delta, \Theta$ as in the previous case. Similarly if $F$ occurs in $\Lambda$. We can therefore assume that $F$ occurs in $\Delta\sharp F$ and in $\Lambda \sharp F$ but it does not occur in $\Delta,\Lambda$.
 Furthermore we can assume that $F$ does not occur in $\Gamma, \Theta$ either, for, otherwise $\Gamma, \Lambda\seq \Delta, \Theta$ can be derived by weakening  ${\cal D} $,  if $F$ occurs in $\Theta$, or ${\cal E}$, if $F$ occurs in $\Gamma$, and then contracting
 the occurrences of $F$ in
 $\Delta \sharp F$ 
with one of the occurrences of $F$
 in $\Theta$ or the occurrences of $F$ in 
$\Lambda \sharp F$ with one of the occurrences of $F$ in $\Gamma$.
Finally if $F$ is atomic it suffices to contract  the occurrence in $F$ in $\Delta\sharp F$  and $\Lambda \sharp F$ into a single one,  and then apply,  possibly after some exchanges, a   cut with  the  atomic cut formula $F$, in order to obtain  the desired derivation.

In the remaining cases  we   proceed, as in Gentzen's original proof of the cut elimination theorem,  by a principal induction on the degree of $F$ and a secondary induction  
 on the sum of the left rank $\rho_l(F,{\cal D})$ of $F$ in ${\cal D}$   and of the right rank $\rho_r(F,{\cal E})$ of $F$ in ${\cal E}$, defined as the largest number of consecutive sequents in a path of ${\cal D}$ (of ${\cal E}$)
starting with the  endsequent, that contain $F$ in the succedent (in the antecedent).

Besides the cases considered in Gentzen's proof, there is also    the possibility that ${\cal D}$ or ${\cal E}$ end  with an atomic equality inference or with an atomic cut.

Case 1  ${\cal D}$ ends, say,  with an atomic  $=_1$-inference. Since $F$ is not atomic,  $F$ is not active in such an inference and   ${\cal D}$ can be represented as:

\[
\ba{c}
{\cal D}_0\\
\Gamma'\seq \Delta'\sharp F, A\{v/r\}\\
\cline{1-1}
\Gamma', r=s\seq \Delta'\sharp F,  A\{v/s\}
\ea
\]
where $\Gamma', r=s$ coincides with $\Gamma$   and  $\Delta', A\{v/s\}$  coincides with $\Delta$. 
Since $\rho_l(F, {\cal D}_0)	<  \rho_l(F, {\cal D})$, by induction hypothesis we have a derivation whose equality and cut-inferences are atomic  of $\Gamma', \Lambda \seq \Delta', A\{v/r\}, \Theta$,
from  which the desired derivation of $\Gamma,\Lambda \seq \Delta, \Theta$ can be obtained by applying the same $=_1$-inference.

Case 2  ${\cal D}$ ends with an atomic cut.  Then ${\cal D}$ can be represented  as:

\[.
\ba{ccc}

{\cal D}_0&~~~~&{\cal D}_1\\
\Gamma_1 \seq \Delta_1\sharp F, A&&\Gamma_2,A\seq \Delta_2\sharp F\\
\cline{1-3}
\multicolumn{3}{c}{\Gamma_1,\Gamma_2 \seq \Delta \sharp F}
\ea
\]where $\Delta$ coincides with $\Delta_1,\Delta_2$, so that $F$ does not occur in  $\Delta-1$ nor in $\Delta_2$.
Since $\rho_l(F, {\cal D}_0) < \rho_l(F, {\cal D})$  
and  $\rho_l(F, {\cal D}_1) < \rho_l(F, {\cal D})$, 
by induction hypothesis applied to ${\cal D}_0$ and ${\cal E}$ and to 
${\cal D}_1$ and ${\cal E}$ there are derivations whose equality and cut-inferences are atomic of  
 $\Gamma_1, \Lambda \seq \Delta_1, A, \Theta$ and 
$\Gamma_2, A, \Lambda \seq \Delta_2, \Theta$, 
to which   it suffices to apply a cut with atomic cut formula $A$ ad then some exchanges and contraction  to have the desired derivation of 
$\Gamma_1,\Gamma_2,\Lambda \seq \Delta_1,\Delta_2,\Theta$. 

The cases in  which it is ${\cal E}$ to end with an atomic equality or  a cut-inference are entirely analogous.
$\Box$.

\ 

From Proposition \ref{atomic} and Proposition \ref{atomic-cut} it follows immediately the following:

\begin{proposition} \label{atomic-atomic}
Every derivation in $LJ^=$ or $LK^=$ can be transformed into a derivation of its endsequent, whose equality and cut-inferences are atomic.
\end{proposition}

{\bf Remark } For the proof of Proposition \ref{atomic-cut} it is crucial that the equality rules transform  atomic formulae only.
For example  in case $\rho_l(F, {\cal D})=1$ and $\rho_r(F,{\cal E})= 1$, if $F$ had the form $F^\circ\{v/s\}$,
with $F^\circ$ non atomic,  ${\cal D}$ ended  with an equality inference transforming $F^\circ\{v/r\}$ into $ F^\circ\{v/s\}$, and ${\cal E}$ by a logical inference introducing $F^\circ\{v/s\}$ in the antecedent, then there would be no way of applying the induction hypothesis. 

\

{\bf Note} Concerning the use of $\Gamma\sharp F$, we note that 
when $\Delta\sharp F$ and $\Lambda \sharp F$ take the form $\Delta, F$ and $\Gamma, F$,  from Proposition \ref {atomic-cut}, we obtain  directly that the  derivations
having only atomic equality and cut-inferences are closed under the application of the cut rule. That is a slight simplification with respect to the use of Gentzen's mix rule that eliminates all the occurrences of $F$, so that the use of additional weakenings and exchanges may be necessary to derive the conclusion of a  cut-inference.

\begin{proposition}\label{sepeq}
If  $~\Gamma\seq \Delta \sharp A\{v/r\}$  has a separated derivation in $LJ^=$ or $LK^=$, then also
$\Gamma, r=s \seq \Delta, A\{v/s\}$  and $\Gamma, s=r \seq \Delta, A\{v/s\}$ have separated derivations in the same system.

\end{proposition}

{\bf Proof} Let ${\cal D}$ be a separated derivation of $\Gamma\seq \Delta\sharp A\{v/r\}$. We proceed  by induction on the height $h({\cal D})$ of ${\cal D}$. In the base case ${\cal D}$ reduces to an axiom  and 
it suffices to apply an $=_1$ or an $=_2$-inference to the axiom itself. If $h({\cal D})>0$ we have the following cases.

Case 1.   ${\cal D}$ ends with  a cut or an equality-inference. In this case  ${\cal D}$ doesn't contain any logical inference. If $\Delta=\Delta \sharp A\{v/r\}$, then it suffices to weaken the endsequent of ${\cal D}$. Otherwise we can contract all the occurrences of $A\{v/r\}$ not belonging to $\Delta$ into a single one and then apply an $=_1$ or $=_2$-inference.
 
 Case 2  ${\cal D}$ ends with a weak structural inference. If such an inference  involves one of the occurrences of $A\{v/r\}$ in $\Delta\sharp A\{v/r\}$ not belonging to $\Delta$,
then the desired derivation is provided directly by the induction hypothesis. Otherwise the latter is  obtained by applying the induction hypothesis and then the same weak structural rule.
  
 Case 3  ${\cal D}$ ends with a logical rule.  $A\{v/r\}$, being  atomic, cannot be the principal  formula of the inference, and the conclusion is a straightforward consequence of the induction hypothesis. For example if ${\cal D}$ has the form:

\[
\ba{ccc}
{\cal D}_0&~~~~&{\cal D}_1\\
\Gamma_0\seq \Delta_0\sharp A\{v/r\}, F &&\Gamma_1, G\seq \Delta_1\sharp A\{v/r\}\\
\cline{1-3}
\multicolumn{3}{c}{ \Gamma_0,\Gamma_1, F\imp G\seq \Delta  \sharp A\{v/r\} }
\ea
\]
where $\Gamma$ coincides with $\Gamma_0,\Gamma_1, F\imp G$,and $\Delta$ with $\Delta_0,\Delta_1$,  by induction hypothesis we have separated derivation ${\cal D}'_0$ and 
${\cal D}_1'$ of
$\Gamma_0, r=s\seq \Delta_0, A\{v/s\}, F$  and $\Gamma_1, G, r=s \seq A\{v/s\}$. 
Then 

\[
\ba{ccc}
{\cal D}'_0&~~~~&{\cal D}'_1\\
\Gamma_0, r=s\seq\Delta_0, A\{v/s\}, F &&\Gamma_1, G, r=s \seq \Delta_1, A\{v/s\}\\
\cline{1-3}
\multicolumn{3}{c}{\Gamma_0,\Gamma_1, r=s, r=s, F\imp G\seq \Delta_0, A\{v/s\}, \Delta_1,  A\{v/s\}}\\
\cline{1-3}
\multicolumn{3}{c}{\Gamma_0,\Gamma_1, r=s, F\imp G\seq \Delta_0,  \Delta_1,  A\{v/s\}}\\
\ea
\]
is  a  separated derivation of  $\Gamma, r=s \seq \Delta, A\{v/s\}$.  $\Box$

\begin{proposition}\label{sepcut}

If $\Gamma \seq \Delta\sharp A$ and $\Lambda \sharp A \seq \Theta$ have separated derivations in $LJ^=$ or $LK^=$,  then also $\Gamma, \Lambda \seq \Delta, \Theta$ has
 a separated derivation in the same system.
\end{proposition}

{\bf Proof}. Let ${\cal D}$ and ${\cal E}$ be separated  derivations of  $\Gamma \seq \Delta\sharp A$ and $\Lambda \sharp A \seq \Theta$ respectively. If $\Delta\sharp A=\Delta$ or $\Lambda \sharp A=\Lambda$ the desired derivation can be obtained by weakening the conclusion of ${\cal D}$ or of ${\cal E}$. If both ${\cal D}$ and ${\cal E}$ end with an equality-inference of with a cut, then ${\cal D}$ and ${\cal E}$, being separated, do not contain any logical inference.
Then  it suffices to contract all the occurrences of $A$ in $\Delta\sharp A$ not occurring in $\Delta$ and, similarly,  all those occurring in $\Lambda \sharp A$ but not in $\Lambda$, into a  single one, and  apply an  atomic cut on $A$.
If ${\cal D}$ or ${\cal E}$, say ${\cal D}$, ends with a weak structural inference or with a logical  inference, we proceed by induction on the sum $h({\cal D}) + h({\cal E})$ of the heights of ${\cal D}$ and ${\cal E}$.

Case 1  ${\cal D}$ ends with a weak structural inference. If such an inference  involve one of the occurrences of $A\{v/r\}$ in $\Delta\sharp A\{v/r\}$ not belonging to $\Delta$,
then the desired derivation is provided directly by the induction hypothesis. Otherwise the latter is  obtained by applying the induction hypothesis and then the same weak structural rule.

Case 2 ${\cal D}$ ends with a logical inference. Since $A$ is atomic, $A$ is not the principal formula of such an inference. Then the conclusion follows by a straightforward induction on $h({\cal D}) + h({\cal E})$. For example if ${\cal D}$ is of the form:

\[
\ba{ccc}
{\cal D}_0&~~~~~&{\cal D}_1\\
\Gamma', F\seq \Delta \sharp A&&\Gamma', G\seq \Delta \sharp A\\
\cline{1-3}
\multicolumn{3}{c}{ \Gamma', F\vel G \seq \Delta \sharp A}
\ea
\]
By induction hypothesis applied to ${\cal D}_0$ and ${\cal E}$ and to ${\cal D}_1$ and  ${\cal E}$ we have two separated derivations of $\Gamma', F, \Lambda \seq \Delta, \Theta$ and $\Gamma', G, \Lambda \seq \Delta, \Theta$, from which by a $\vel \seq$-inference we obtain the desired derivation of $\Gamma', F\vel G \seq \Delta, \Theta$.

The cases in which it is ${\cal E}$ to end with a weak structural inference or with a cut are entirely analogous.
$\Box$

\begin{proposition}\label{redsep}
Every derivable sequent in $LJ^=$ or $LK^=$ has a separated derivation in the same system,

\end{proposition}

{\bf Proof} Assume we are given  a non separated derivation  ${\cal D}$ of $\Gamma \seq \Delta$ in $LJ^=$ or $LK^=$. 
By Proposition \ref{atomic-atomic}, ${\cal D}$ can be transformed into a derivation ${\cal D}'$ whose equality and cut-inferences are atomic. Then a straightforward induction on the height of ${\cal D}'$, based on Proposition \ref{sepeq} and Proposition \ref{sepcut} shows that ${\cal D}'$ can be transformed into a separated derivation of $\Gamma \seq \Delta$.
  $\Box$


\

By the previous Proposition \ref{redsep},  to show that the cut rule is eliminable from derivations in $LJ^=$ or $LK^{=}$  it suffices to show that it can be  eliminated from the derivations of the purely equational calculus $EQ$.
Instrumental for that pourpose will  be the following equational calculus $EQ_N$, where  $N$ stands for {\em natural}.

\begin{definition}
$EQ_N$ is the calculus acting on  sequents with one formula in the succedent,  obtained from $EQ$ by replacing the rules $=_1$ and $=_2$ with the rule
$CNG$:

\[
\ba{c}
\Gamma\seq F\{v/r\}~~~~~~~\Lambda \seq r=s\\
\cline{1-1}
\Gamma, \Lambda\seq F\{v/s\}
\ea
\]

\end{definition}

\begin{definition}
$cf.EQ$ and $cf.EQ_N$ denote the systems $EQ$ and $EQ_N$ deprived of the cut rule.

\end{definition}

\begin{proposition}\label{equiv}

EQ and $EQ_N$ are equivalent.
\end{proposition}

{\bf Proof} The following are derivations of $=_1$ and $=_2$  from $CNG$  and of $CNG$ from $=_1$:
\[
\ba{cccl}
\Gamma\seq F\{v/r\}&~~~&r=s\seq r=s&\vbox to 0pt{\hbox{\scriptsize{CNG}}}\\
\cline{1-3}

\multicolumn{3}{c}{\Gamma, r=s\seq F\{v/s\}}

\ea
\]

\[
\ba{cccl}
&&\seq s=s~~~~s=r\seq s=r& \vbox to 0pt{\hbox{\scriptsize{CNG}}}\\
\cline{3-3}
\Gamma\seq F\{v/r\}&&s=r\seq r=s&\vbox to 0pt{\hbox{\scriptsize{CNG}}}\\
\cline{1-3}
\multicolumn{3}{c}{\Gamma, s=r\seq F\{v/s\}}\
\ea
\]

\[
\ba{cccl}
&~~~&\Gamma \seq F\{v/r\}&\vbox to 0pt{\hbox{$=_1$}}\\
\cline{3-3}
\Lambda \seq r=s &&\Gamma, r=s \seq F\{v/s\}\\
\cline{1-3}
\multicolumn{3}{c}{\Gamma, \Lambda \seq F\{v/s\}}
\ea
\]
$\Box$

\subsection{Cut-elimination for $EQ_N$}

\begin{proposition}\label{cutlemma}
If $\Gamma \seq F$ and $\Lambda \sharp F\seq G$ are derivable in $cf. EQ_N$, then also $\Gamma, \Lambda \seq G$ is derivable in $cf.EQ_N$.
\end{proposition}

{\bf Proof} Let ${\cal D}$  and ${\cal E}$ be  derivations   in $cf.EQ$ of $\Gamma \seq F$, and  $\Lambda \sharp F\seq G$ respectively.
We have to show that  there is a derivation ${\cal F}$ in $cf.EQ $ of $\Gamma, \Lambda \seq G$. 

If  $\Lambda \sharp F$  coincides with $\Lambda$, in particular if $\Lambda\sharp F$ is empty, or $F$ occurs in $\Lambda$, then
to obtain ${\cal F}$ it suffices to add to ${\cal E}$ the weakenings, exchanges  and, in the latter case,  contractions needed to obtain $\Gamma, \Lambda \seq G$.
Otherwise we proceed by induction on the height $h({\cal E})$ of ${\cal E}$. If $h({\cal E})=0$, then 
${\cal E}$ reduces to $F\seq F$ and for ${\cal F}$ we can take ${\cal D}$ itself.

If ${\cal E}$ ends with a weak structural inference  that involves (at least) one of the occurrences of $F$ in $\Lambda\sharp F$, that does not occur in $\Lambda$, then the desired derivation ${\cal F}$ is provided directly by the induction hypothesis.
Otherwise it suffices to apply the induction hypothesis and then the last weak structural inference of ${\cal E}$.

If  ${\cal E}$ ends with a $CNG$-inference, then $G$ has the form $H\{v/s\}$ and ${\cal E}$  can be represented as:

\[
\ba{ccc}
{\cal E}_0&~~~&{\cal E}_1\\
\Lambda_0\sharp F \seq H\{v/r\}&&\Lambda_1\sharp F \seq r=s\\
\cline{1-3}
\multicolumn{3}{c}{ \Lambda \sharp F \seq H\{v/s\}}
\ea
\]
By induction hypothesis we have cut-free derivations of $\Gamma, \Lambda_0 \seq H\{v/r\}$ and 
 $\Gamma , \Lambda_1 \seq r=s$, from which ${\cal F}$ is obtained by applying  the same $CNG$-inference
and some exchanges and contractions. $\Box$

\begin{proposition}\label{EQNcutelimination}
If  a sequent  is derivable in $EQ_N$, then it is also derivable in $cf.EQ_N$.
\end{proposition}

{\bf Proof} By the previous Proposition,  applied in the specific case in which $\Lambda\sharp F$ is $\Lambda, F$, it follows that  the   cut rule
is admissible in $cf.EQ_N$
and therefore eliminable from derivations in $EQ_N$. $\Box$

\subsection{Cut elimination for $LJ^=_N$ and $LK^=_N$}

From Proposition \ref{redsep} , 
Proposition \ref{equiv} 
and Proposition \ref{EQNcutelimination} we obtain the  full cut elimination theorem for the calculi $LJ_N^=$ and $LK_N^{=}$, that are obtained by adding to  $LJ$ and $LK$ 
the Reflexivity Axiom $\seq =$   and the rule $CNG$.

\begin{theorem}\label{cutL}
The cut rule is eliminable from derivations in $LJ_N^=$ and in $LK_N^{=}$.
\end{theorem}

\subsection{Admissibility of $CNG$ in $cf.EQ$}

\begin{proposition} \label{CNGadmissibility} The rule $CNG$ is admissible in $cf.EQ$, namely, 
if $\Gamma \seq F\{v/r\}$  and $\Lambda \seq r=s$ are derivable in $cf.EQ$ then also $\Gamma,\Lambda \seq F\{v/s\}$ is derivable in $cf.EQ$.
\end{proposition}

{\bf Proof} Let  ${\cal D}$ and  ${\cal E}$  be derivations in $cf.EQ$ of  $\Gamma \seq F\{v/r\}$ and  $\Lambda \seq r=s$ respectively. 
We have to show that there is a derivation  ${\cal F}$ of  $\Gamma, \Lambda \seq F\{v/s\}$ in $cf.EQ$.
If  $r$ and  $s$ coincide, to obtain  ${\cal F}$ it suffices to  apply to the end-sequent  of  ${\cal D}$   the appropriate weakenings to introduce $\Lambda$ in the antecedent
 of its end-sequent. Otherwise we proceed by induction on the height of  ${\cal E}$, with respect to an arbitrary ${\cal D}$.  In the  base  case ${\cal E}$ reduces to the axiom  $r=s\seq r=s$.
  In that case as  ${\cal F}$ we can take:

\[
\ba{c}
{\cal D}\\
\Gamma \seq F\{v/r\}\\
\cline{1-1}
\Gamma, r=s \seq F\{v/s\}
\ea
\]that uses $=_1\seq$.
If  ${\cal E}$ ends with a structural rule, to obtain  ${\cal F}$  it suffices to apply the induction hypothesis to ${\cal D}$ and to the immediate subderivation
 ${\cal E}_0$ of ${\cal E}$ and then the  last structural rule of ${\cal E}$. 

If  ${\cal E}$ ends with  a $=_1\seq $-inference, namely it is of the form:

\[
\ba{c}
{\cal E}_0\\
\Lambda' \seq r^\circ\{u/p\}=s^\circ\{u/p\}\\
\cline{1-1}
\Lambda', p=q \seq r^\circ\{u/q\}=s^\circ\{u/q\}
\ea
\]
so that  $r$ and  $s$ are  $r^\circ\{u/q\}$ and  $s^\circ \{u/q\}$ respectively, and  $\Lambda $ is  $\Lambda', p=q$,  let ${\cal D}'$ be the following derivation:

\[
\ba{c}
{\cal D}\\
\Gamma \seq F\{v/r^\circ\{u/q\}\}\\
\cline{1-1}
\Gamma, p=q  \seq F\{v/r^\circ\{u/p\}\}
\ea
\]
which uses $=_2$. By induction hypothesis applied to 
${\cal D}'$ and  ${\cal E}_0$ there is a derivation  ${\cal F}_0$ of  $\Gamma, p=q, \Lambda' \seq F\{v/s^\circ \{ u/p\}\}$.
As  ${\cal F}$ we can then take the following derivation:

\[
\ba{c}
{\cal F}_0\\
\Gamma, p=q, \Lambda' \seq F\{v/s^\circ \{ u/p\}\}\\
\cline{1-1}
\Gamma, p=q, \Lambda', p=q\seq F\{v/s^\circ \{ u/q\}\}\\
\cline{1-1}
\Gamma,  \Lambda', p=q \seq F\{v/s^\circ \{ u/q\}\}\\
\ea
\]
which uses $=_1\seq$ and a contraction.

Finally if  ${\cal E}$ ends with a $=_2$-inference, namely it is of the form:

\[
\ba{c}
{\cal E}_0\\
\Lambda' \seq r^\circ\{u/p\}=s^\circ\{u/p\}\}\\
\cline{1-1}
\Lambda', q=p \seq   r^\circ\{u/q\}=s^\circ\{u/q\}\}
\ea
\]

We let   ${\cal D}'$ be :

\[
\ba{c}
{\cal D}\\
\Gamma \seq F\{v/r^\circ\{u/q\}\}\\
\cline{1-1}
\Gamma,q=p \seq F\{v/r^\circ\{u/p\}\}
\ea
\]
 which uses $=_1\seq$. By induction hypothesis applied to  ${\cal D}'$  e ${\cal E}_0$ we obtain a derivation  ${\cal F}_0$ of   $\Gamma, q=p, \Lambda' \seq F\{v/s^\circ \{ u/p\}\}$.
 Then, as  ${\cal F}$ we take the following derivation:
 
 \[
 \ba{c}
 {\cal F}_0\\
 \Gamma, q=p, \Lambda' \seq F\{v/s^\circ \{ u/p\}\}\\
 \cline{1-1}
  \Gamma, q=p, \Lambda' , q=p \seq F\{v/s^\circ \{ u/q\}\}\\
  \cline{1-1}
  \Gamma, \Lambda' , q=p \seq F\{v/s^\circ \{ u/q\}\}
 \ea
 \]
  which uses  $=_2$ and a contraction. $\Box$

\subsection{Cut elimination for $EQ$}
 
 \begin{theorem} \label{CutEQ}
Cut elimination for $EQ$.

If $\Gamma \seq F$ is derivable in $EQ$, then it is derivable also  in $cf.EQ$
 \end{theorem}
 
 {\bf Proof} 
 By Proposition \ref{equiv} a derivation  ${\cal D}$ of  $\Gamma \seq F$ in $EQ$
 can be transformed into a derivation  ${\cal D}'$  in $EQ_N$ of  $\Gamma \seq F$.
 By the eliminability of the cut-rule in $EQ_N$,
${\cal D}'$ can be transformed into a derivation ${\cal D}''$ in $cf.EQ_N$ of $\Gamma\seq F$. Finally by the admissibility of $CNG$ in $cf.EQ$,  ${\cal D}''$ can be transformed into a derivation  in $cf.EQ$ of $\Gamma \seq F$. $\Box$

\subsection{Cut elimination for $LJ^=$ and $LK^{=}$}

From Proposition \ref{redsep} and Theorem\ref{CutEQ} we obtain the  full cut elimination theorem for $LJ^=$ and $LK^{=}$.

\begin{theorem}\label{cutL}
The cut rule is eliminable from derivations in $LJ^=$ and in $LK^{=}$.
\end{theorem}

\subsection{Cut elimination for $LJ^{(1)=}$ and $LK^{(1)=}$}

Since,  the rules $=_1$ and $=_2$ are derivable in $LJ^{(1)=}$,  without using the cut rule, from cut elimination for $LJ^=$
it follows immediately that cut elimination holds also for $LJ^{(1)=}$ and $LK^{(1)=}$:

\begin{theorem}\label{cutL1}
The cut rule is eliminable from derivations in $LJ^{(1)=}$ and in $LK^{(1)=}$.
\end{theorem}

\subsection{Admissibility of $=_1^l$ and $=_2^l$ in $cf.EQ$}

Since $=_1^l$ and $=_2^l$, namely:

\[
\ba{clccl}
\Gamma, F\{v/r\}\seq \Delta&\vbox to 0pt{\hbox{$=_1^l$ }}&~~~\mbox{and}~~~&\Gamma, F\{v/r\}\seq \Delta&\vbox to 0pt{\hbox{$=_2^l$ }}\\
\cline{1-1}\cline{4-4}
\Gamma, F\{v/s\}, r=s\seq \Delta&&&\Gamma, F\{v/s\}, s=r\seq \Delta
\ea
\]
 are derivable in $EQ$ and the cut rule is eliminable from derivations in $EQ$ we immediately have the following:

\begin{proposition}\label{l-admissibility}

The rules $=_1^l$ and $=_2^l$ are admissible in $cf.EQ$.

\end{proposition}\label{L-admissibility}

\begin{definition}
Let $EQ_{1}$ be   obtained from $EQ$ by replacing $=_2$ by $=_1^l$ and   $EQ_{2}$  be    obtained from $EQ$ by replacing $=_1$ by $=_2^l$ .
 $cf.EQ_{1}$ and $cf.EQ_{2}$ denote  $EQ_{1}$ and $EQ_{2}$  deprived of the cut rule.
\end{definition}

\subsection{Admissibility of $=_2$  in $EQ_{1}$ and of $=_1$ in $EQ_{2}$}

{\bf Notation} In the following $E\equiv E'$ will denote syntactic  equality between the terms or formulae  that are denoted by $E$ and $E'$.

\ 

As already noted in \cite{L68} we have the following:

\begin{lemma}\label{singleton}
The equality rules $=_1$ and $=_2$ as well as $=_1^l$ and $=_2^l$ are derivable by means of the contraction rule from their {\em singleton}
version, obtained by requiring that $v$ has exactly one occurrence in the changing formula.
\end{lemma}

{\bf Proof} It suffices to deal with $=_1$, the other cases being entirely similar.

Given $F$ with $n$ occurrence  of $v$, with $n>1$,  let $F'$ be obtained from $F$ by replacing all the  occurrences of $v$ by $n$ new (to $F$, $r$ and $s$) distinct 
variables $v_1,\ldots, v_n$, so that $F\{v/r\}\equiv F'\{v_1/r,\ldots, v_n/r\}$ and $F\{v/s\}\equiv F'\{v_1/s,\ldots, v_n/s\}$.
\[
\ba{c}
\Gamma\seq F'\{v_1/r, \ldots, v_{n-1}/r, v_n/r\}\\
\cline{1-1}
\Gamma, r=s \seq F'\{v_1/r,\ldots, v_{n-1}/r, v_n/s\}
\ea
\]
is a correct application of the singleton version of $=_1$, since $F'\{v_1/r,\ldots, v_{n-1}/r_{n-1}, v_n/r_n\}\equiv F'\{v_1/r,\ldots, v_{n-1}/r\}\{v_n/r\}$
and $  F'\{v_1/r,\ldots, v_{n-1}/r\}\{ v_n/s \} \equiv F'\{v_1/r,\ldots, v_{n-1}/r, v_n/s\}$. Similarly,
since  $F'\{v_1/r,\ldots, v_{n-2}/r, v_{n-1}/r, v_n/s\}\equiv F'\{v_1/r,\ldots, v_{n-2}/r,  v_n/s\}\{v_{n-1}/r\}$ 
and $F'\{v_1/r,\ldots, v_{n-2}/r,  v_n/s\}\{v_{n-1}/s\}\equiv F'\{v_1/r,\ldots, v_{n-2}/r, v_{n-1}/s, v_n/s\}$
 the following  it is a correct application of $=_1$:
\[
\ba{c}
\Gamma, r=s \seq F'\{v_1/r,\ldots, v_{n-2}/r, v_{n-1}/r, v_n/s\}\\
\cline{1-1}
\Gamma, r=s,r=s \seq  F'\{v_1/r,\ldots, v_{n-2}/r, v_{n-1}/s,  v_n/s\}\
\ea
\]
Proceeding in that way, with $n$ applications of the singleton  $=_1$-rule we obtain a derivation from $\Gamma\seq F\{v/r\}$ of
$\Gamma, r=s,\ldots, r=s\seq F\{v/s\}$, from which the desired derivation of $\Gamma, r=s \seq F\{v/s\}$ can be obtained by  $n-1$ applications of the contraction rule. $\Box$

\begin{definition}

$cf.EQ^1_{1}$  and $cf.EQ^1_{2}$,  are obtained from  $cf.EQ_{1}$  and $cf.EQ_{2}$
by replacing the equality rules by their singleton version.
\end{definition}

\begin{proposition}\label{2admissibility}
$=_2 $ is admissible in $cf.EQ_{1}$ and $=_1$ is admissible in $cf.EQ_{2}$.

\end{proposition}
{\bf Proof} By the previous Lemma \ref{singleton} it suffices to prove that the singleton versions of $=_2$ and $=_1$ are admissible in the systems $cf.EQ^1_{1}$  and $cf.EQ^1_{2}$, namely that:

$a)$  if $\Gamma \seq F\{v/r\}$ is derivable in $cf.EQ^1_{1}$,   then also $\Gamma, s=r\seq F\{v/s\}$ is derivable in $cf.EQ^1_{1}$,
and 

$b)$   if $\Gamma \seq F\{v/r\}$ is derivable in $cf.EQ^1_{2}$,   then also $\Gamma, r=s\seq F\{v/s\}$ is derivable in $cf.EQ^1_{2}$.

As for $a)$, let ${\cal D}$ be a derivation in $cf.EQ^1_{1}$ of $\Gamma \seq F\{v/r\}$. We proceed by induction  on the height $h({\cal D})$
of ${\cal D}$  to show that in $cf.EQ^1_{1}$ there is a derivation ${\cal D}'$ of $\Gamma, s=r\seq F\{v/s\}$.
If $h({\cal D})=0$ then ${\cal D}$ reduces to $F\{v/r\}\seq F\{v/r\}$ or to $\seq t_0=t_1\{v/r\}$, with  $t_0\equiv  t_1\{v/r\}$ or to $\seq t_0\{v/r\} = t_1$, with $t_0\{v/r\} \equiv t_1$
In the former case as ${\cal D}'$ we can take:
\[
\ba{cl}
F\{v/s\}\seq F\{v/s\}&\vbox to 0pt{\hbox{$=_1^l $}}\\
\cline{1-1}
F\{v/r\},s=r\seq F\{v/s\}
\ea
\]
If ${\cal D}$ reduces to $\seq t_0=t_1\{v/r\}$, with  $t_0\equiv  t_1\{v/r\}$ as ${\cal D}'$ we can take:

\[
\ba{cl}
\seq t_1\{v/s\}=t_1\{v/s\}&\vbox to 0pt{\hbox{$=_1 $}}\\
\cline{1-1}
s=r\seq t_0 = t_1\{v/s\}
\ea
\]which is   correct  since $t_0\equiv t_1\{v/r\}$.
The case in which  ${\cal D}$ reduces to $\seq t_0\{v/r\}=t_1$, with $t_0\{v/r\}\equiv t_1$, is entirely similar.



If $h({\cal  D})>0$, and ${\cal D}$ ends  with a  structural rule  the conclusion is a straightforward
consequence of the induction hypothesis.
If ${\cal D}$ ends with a $=_1$-inference, then  we distinguish the following three subcases.

Case 1. ${\cal D}$ is of the form:
\[
\ba{cl}
{\cal D}_0\\
\Gamma'\seq F^\circ\{u/p, v/r\}&\vbox to 0pt{\hbox{$=_1 $}}\\
\cline{1-1}
\Gamma', p=q \seq F^\circ \{u/q, v/r\}\
\ea
\]
with  $F\equiv F^\circ\{u/q\}$ and the unique occurrence of $v$ in $F$ does not occur in $q$ (and $\Gamma$ coincides with $ \Gamma', p=q$).

By induction hypothesis we have a derivation ${\cal D}_0'$ in $EQ^1_{1}$ of $\Gamma', s=r \seq F^\circ\{u/p, v/s\}$.
As ${\cal D}'$ we can then take:

\[
\ba{cl}
{\cal D}'_0\\
\Gamma', s=r \seq F^\circ \{u/p, v/s\}&\vbox to 0pt{\hbox{$=_1 $}}\\
\cline{1-1}
\Gamma', p=q, s=r \seq F^\circ \{u/q, v/s\}\
\ea
\]
Case 2. 
 ${\cal D}$ is of the form:
\[
\ba{cl}
{\cal D}_0\\
\Gamma'\seq F^\circ\{u/p\}&\vbox to 0pt{\hbox{$=_1 $}}\\
\cline{1-1}
\Gamma', p=q\{v/r\} \seq F^\circ\{u/q\{v/r\}\}\
\ea
\]
with $F\equiv F^\circ\{u/q\}$ and the unique occurrence of $v$ in $F$ occurs in $q$.

As ${\cal D}'$ we can then take:

\[
\ba{cl}
{\cal D}_0&\\
\Gamma'\seq F^\circ \{u/p\}&\vbox to 0pt{\hbox{$=_1 $}}\\\
\overline{\Gamma', p=q\{v/s\}\seq F^\circ\{u/q\{v/s\}\}}&\vbox to 0pt{\hbox{$=_1^l $}}\\
\overline{\Gamma', p=q\{v/r\}, s=r \seq F^\circ \{u/q\{v/s\}\}}&
\ea
\]

Case 3. ${\cal D}$ is of the form:
\[
\ba{cl}
{\cal D}_0\\
\Gamma'\seq F\{v/r^\circ\{u/p\} \}&\vbox to 0pt{\hbox{$=_1 $}}\\
\cline{1-1}
\Gamma', p=q \seq F\{v/ r^\circ\{u/q\}\}\
\ea
\]and  $r\equiv r^\circ\{u/q\}$.
By induction hypothesis we have a derivation ${\cal D}_0'$ in $EQ^1_1$ of

 $\Gamma', s=r^\circ\{u/p\} \seq F \{v/s\}$.
As ${\cal D}'$ we can take:

\[
\ba{cl}
{\cal D}_0'&\\
\Gamma', s=r^\circ\{u/p\} \seq F\{v/s\}&\vbox to 0pt{\hbox{$=_1^l$}}\\
\cline{1-1}
\Gamma', s=r^\circ\{u/q\}, p=q \seq F\{v/s\}
\ea
\]
If ${\cal D}$ ends with a $=_1^l$-inference, then ${\cal D}$ has the  form:

\[
\ba{cl}
{\cal D}_0\\
\Gamma', G\{u/p\}\seq F\{v/r\}&\vbox to 0pt{\hbox{$=_1^l $}}\\
\cline{1-1}
\Gamma', G\{u/q\}, p=q \seq F\{v/r\}
\ea
\]
By induction  hypothesis we have a derivation  ${\cal D}_0'$ in $EQ_1^1$ of
$\Gamma', G\{u/p\}, s=r \seq F\{v/s\}$. As ${\cal D}'$ we can take:

\[
\ba{cl}
{\cal D}'_0&\\
\Gamma', G\{u/p\}, s=r\seq F\{v/s\}& \vbox to 0pt{\hbox{$=_1^l$}}\\
\cline{1-1}
\Gamma', G\{u/q\}, p=q, s=r \seq F\{v/s\}
\ea
\]The proof of $b)$ is entirely similar.$\Box$

\begin{theorem}\label{cutelimination$EQ_1$}
Cut elimination holds for $EQ_1$ and $EQ_2$.
\end{theorem}

{\bf Proof} Any derivation  ${\cal D}$ in $EQ_1$ can be transformed (by using the cut rule) into a derivation ${\cal D}'$ in $EQ$ of the same end sequent.
By the cut elimination theorem for $EQ$, ${\cal D}'$ can be transformed into a cut-free derivation ${\cal D}"$ in $EQ$. Since $=_2$ is admissible in $cf.EQ_1$,
the applications of the $=_2$-rule  in ${\cal D}"$ can be replaced by applications of $=_1$ and $=_1^l$, thus obtaining the desired cut free derivation
in $EQ_1$  of the end sequent of ${\cal D}$. Thanks to the admissibility in $EQ_2$ of $=_1$, the same argument shows that cut elimination holds for $EQ_2$ as well. $\Box$

\ 

{\bf Note} Since $=_2^l$ is derivable in $EQ_1$, from the admissibility of the cut rule in $cf.EQ_1$
it follows that $=_2^l$ is also admissible in $cf. EQ_1$. Similarly also $=_1^l$ is admissible in $EQ_2$.

\begin{definition}
For $i=1,2$, $LJ^=_{i}$ and $LK^=_{i}$ denote the systems obtained by adding $=_i$ and $=_i^l$ to $LJ$ and $LK$ respectively. 

\end{definition}

As an immediate consequence of Theorem \ref{cutelimination$EQ_1$}, we have the following
\begin{theorem}
Cut elimination holds for $LJ^=_1$, $LJ^=_2$, $LK^=_1$ and $LK^=_2$.
\end{theorem}

$EQ$, $EQ_1$ and $EQ_2$ are the only systems satisfying cut elimination  that can  be obtained by adding  to the structural rules the Reflexivity Axiom  and  two equality rules chosen among $=_1$, $=_2$, $=_1^l$ and  $=_2^l$.

 For example

$a=c, b=c \seq a=b$ has the following cut-free derivations:
\[
\ba{clccl}
a=c\seq a=c&\vbox to 0pt{\hbox{$=_2$}}&~~~\mbox{and}~~~~&a=b\seq a=b&\vbox to 0pt{\hbox{$=_1^l$}}\\
\cline{1-1}\cline{4-4}
a=c, b=c\seq  a=b&&&a=c,b=c\seq a=b
\ea
\]
but it has no cut-free derivation, if $a$, $b$ and $c$ are distinct and  only the use of $=_1$ and $=_2^l$ is allowed.
More generally no sequent of the form 
$
~*)~\Gamma  \seq a=b
$,
 where the  formulae in $\Gamma$ are among $c=c$, $a=c$ and $b=c$,
  can have a cut free derivation using only $=_1$ and $=_2^l$ .
In fact, $*)$ is not the conclusion of  a non trivial $=_1$-inference, since  $c$ occurs in the right-hand side of all the possible operating equalities, so that it would occur in the  succedent of  the conclusion of any such inference.
If it is the conclusion of a $=_2^l$-inference, with operating equality  $a=c$, the transformed formula must be necessarily another occurrence of $a=c$, obtained by replacing with $a$ the first  occurrence of $c$ in the changing formula $c=c$,  to be found in the antecedent of the premiss. The same holds if the operating equality is 
$b=c$. Thus the premiss of the inference is still a sequent of the form $*)$. Obviously that is the case if $*)$ is the conclusion of a  weakening, exchange or contraction.
Hence if  $*)$ is the conclusion of an inference  different from a cut, then also the premiss of the inference has the form $*)$.
 Assuming that $a$, $b$ and $c$ are distinct,
no axiom has the form $*)$. Thus there are no derivations of height zero of sequents of that   form.  Furthermore, by the above discussion,  if  there are no derivations of height $n$ of sequents of 
the form $*)$,  
then there are no   derivations of height $n+1$  of sequents of 
that same form either.
By induction on $n$ we conclude that there are no derivations at all of sequents
 of the  form $*)$.
In particular, if $a$, $b$ and $c$ are distinct,
 $a=c,b=c \seq a=b$ has no cut-free derivation 
 using only $=_1$ and $=_2^l$ .

Similarly $c=b,c=a\seq a=b$ has the cut-free derivations:

\[
\ba{clccl}
c=b\seq c=b&\vbox to 0pt{\hbox{$=_1$}}&~~~\mbox{and}~~~~&a=b\seq a=b&\vbox to 0pt{\hbox{$=_2^l$}}\\
\cline{1-1}\cline{4-4}
c=b,c=a\seq  a=b&&&c=b,c=a\seq a=b
\ea
\]
but it has no cut-free derivation, if $a$, $b$ and $c$ are distinct and only the use of $=_2$ and $=_1^l$ is allowed.

 Finally  $a=b\seq f(a)=f(b)$ has the cut-free derivations:

\[
\ba{clccl}
\seq f(a)=f(a)&\vbox to 0pt{\hbox{$=_1$}}&~~~\mbox{and}~~~~&\seq f(b)=f(b)&\vbox to 0pt{\hbox{$=_2$}}\\
\cline{1-1}\cline{4-4}
a=b\seq  f(a)=f(b)&&&a=b\seq f(a)=f(b)
\ea
\]
but, if $a$ and $b$ are distinct,  it has no cut-free derivation using only  $=_1^l $ and $=_2^l$. In fact,  if $a$ and $b$ are distinct,
no sequent of the form 
$
\Gamma  \seq f(a)=f(b)
$,
 where the  formulae in $\Gamma$ are among $a=a$, $ b=b$, $a=b$ and $b=a$,
  can have a cut free derivation using only $=_1^l$ and $=_2^l$.

 Clearly that remains the case even if we add 
  the {\em left symmetry rule}, that leads  from $\Gamma, r=s\seq \Delta$,  to   $\Gamma, s=r\seq \Delta$.
Concerning such a rule,  we note also  that it  has the following cut-free derivation based on $=_1^l$, $=_2^l$ and the contraction rule:
  \[
  \ba{cl}
  \Gamma, r=s\seq \Delta&\vbox to 0pt{\hbox{$=_1^l$}}\\
  \cline{1-1}
  \Gamma, r=r, s=r \seq \Delta&\vbox to 0pt{\hbox{$=_2^l$}}\\
  \cline{1-1}
  \Gamma, s=r, s=r, s=r \seq \Delta&\\
  \cline{1-1}
  \Gamma,s=r\seq \Delta&
  \ea
  \]
while it  is not even admissible in the cut-free system with $=_1$ and $=_2^l$.   For,  otherwise,  also $=_1^l$ would be admissible and then  cut elimination would hold, which we have shown not to be the case. Similarly the left symmetry rule is not admissible in the cut free system with $=_1^l$ and $=_2$. On the contrary, since it is derivable in $EQ_1$, $EQ_2$ and $EQ$ (by means of the cut rule)  as shown by the derivations:

\[
\ba{clc}
\seq s=s&\vbox to 0pt{\hbox{$=_1$}}~~&\\
\cline{1-1}
s=r\seq r=s&&\Gamma, r=s\seq \Delta\\
\cline{1-3}
\multicolumn{3}{c}{\Gamma, s=r\seq \Delta}
\ea
\]

\[
\ba{clc}
\seq r=r&\vbox to 0pt{\hbox{$=_2$}}~~&\\
\cline{1-1}
s=r\seq r=s&&\Gamma, r=s\seq \Delta\\
\cline{1-3}
\multicolumn{3}{c}{\Gamma, s=r\seq \Delta}
\ea
\] 
it is admissible in the cut-free part of any of these systems (a  fact that can also be easily proved directly   by induction on the height of derivations).
On the other hand the  left symmetry rule is not derivable in any of $cf.EQ$, $cf.EQ_1$ and $cf.EQ_2$. For $cf.EQ$ that is obvious since
 $=_1$ and $=_2$ add formulae in the antecedent and modify only the formula in the succedent of a sequent. As for $cf.EQ_1$ ($cf.EQ_2$)  it suffices to note that 
all the sequents
in a derivation
that starts  with a sequent  containing $a=b$  in the antecedent, 
  must contain an equality of the form $a=t$ ($t=b$) in the antecedent. As a consequence, for example, there cannot be any derivation in $cf.EQ_1$ or $cf.EQ_2$
of $b=a\seq c=d$ from $a=b\seq c=d$, with $a,b,c$ and $d$ distinct.

\ 

Since any of the four equality rules is derivable from any  other,  from the above discussion concerning the failure of cut elimination, it follows that if only one of them is added to the logical and  reflexivity axioms  and  the structural and logical rules, then  
the system that is obtained is adequate for first-order logic with equality,  but it does not satisfy  cut elimination. 
On the other hand if at least three  of them are added, then cut elimination holds.
More precisely  we have established the following result.

\begin{theorem} \label{cutelimination}
Any extension of $LJ$ or $LK$ obtained by  adding the Reflexivity Axiom $\seq =$ and
some of the rules $=_1, =_2, =_1^l$ and $=_2^l$ is adequate for intuitionistic or classical first order logic with equality, but it 
satisfies the cut elimination theorem if and only if it contains (at least) either both $=_1 $and $=_2$, 
or both $=_1$ and $=_1^l$ or both $=_2 $ and $=_2^l$.
\end{theorem}

\subsection{The Semishortening Property}

Letting $LJ^=_{12}$ be the union of $LJ^=_1$ and $LJ^=_2$ and, similarly, $LK^=_{12}$ be the union of $LK^=_1$ and $LK^=_2$, 
by the previous Theorem, cut elimination holds for both $LJ^=_{12}$ and $LK^=_{12}$. 
On the  ground of the exchange  and contraction rules only, $LK^=_{12}$ is equivalent to the system $G^e$
 in \cite{L68}, which generalizes the rules $=_1$ and $=_1^l$ 
  by permitting  the substitution of $r$ by $s$ in more than one formula and merges them
into a single rule of the form:
\[
\ba{cl}
\Gamma\{v/r\}\seq \Delta\{v/r\}&\\
\cline{1-1}
\Gamma\{v/s\}, r=s \seq \Delta\{v/r\}
\ea
\]
and, similarly, generalizes and merges the rules  $=_2$ and $=_2^l$ into:

\[
\ba{cl}
\Gamma\{v/r\}\seq \Delta\{v/r\}&\\
\cline{1-1}
\Gamma\{v/s\}, s=r \seq \Delta\{v/r\}
\ea
\]
Thus, as an immediate consequence of Theorem \ref{cutelimination},  cut elimination holds for $G^e$.
Actually \cite{L68} deals only  with cut-free derivations in $G^e$ and shows that they can  be transformed into cut-free derivations that do not contain terms that are {\em longer} than those occurring in the end-sequent, under various notion of {\em length} of a term. Clearly a cut-free derivation may contain terms longer than those occurring in the endsequent only if
it contains some equality inference that is {\em lengthening} in the sense that the term $r$ in the premiss is longer that the term 
$s$ by which  it is replaced in the conclusion of the inference. If we let $s\prec r$ to mean that $r$ is longer that $s$, the result in 
$\cite{L68}$ applies to all the binary relation $\prec$ on terms that are strict partial orders congruent  with respect to substitution, namely $r\prec s$ entails $t\{v/r\} \prec t\{v/s\}$, for any term $r$, $s$ and $t$. \cite{O69} states  that  it suffices to require that $\prec$ be antireflexive.
We will base our definitions on such a  weaker requirement and prove  a stronger result, namely that any derivation in $LJ^=_{12}$ or $LK^=_{12}$ can be transformed into one whose equality inferences are all non lengthnening, 
while those of the form $=_1$ and $=_1^l$, or, alternatively, those of the form $=_2$ and $=_2^l$, are actually shortening, namely
satisfy the stronger condition $r\prec s$.  It will suffice to deal with the former case, since the latter is completely symmetric.
In the following $\prec$ will be a fixed, but arbitrary binary antireflexive relation on terms, namely for any term $r$ and $s$, 
 $r\prec s$ entails $s\not \prec r$.

\begin{definition} An application of an $=_1$-inference or of an $=_1^l$-inference with operating equality $r=s$ (or of an application of an $=_2$-inference or of an $=_2^l$-inference with operating equality $s=r$)
 is said {\em nonlengthening} if  $s\not\prec r$ and {\em shortening} if $r\prec s$.
A derivation is said to be  nonlengthening if all its equality inferences are  nonlengthening
and semishortening if it is nonlengthening and, furthermore, all its 
 $=_1$ and $=_1^l$-inferences are shortening.
\end{definition}

\begin{proposition} \label{gamma}
The equality rules $=_1$ and $=_2$ are admissible in   $cf.EQ_{12}$ restricted to semishortening derivations.
More precisely, there are two effective operations  ${\cal G}_1$ and ${\cal G}_2$ such that:
\bi
\item[a)]  if ${\cal D}$ is a semishortening derivation in $cf.EQ_{12}$ of $\Gamma \seq F\{v/r\}$, then for any term $s$,  ${\cal G}_1({\cal D}, r,s)$ is a semishortening  derivation in $cf.EQ_{12}$ of 

$\Gamma, r=s\seq F\{v/s\}$ and
\item[b)] if ${\cal D}$ is a semishortening derivation in $cf.EQ_{12}$ of $\Gamma \seq F\{v/r\}$,  then for any term $s$,  ${\cal G}_2({\cal D},r,  s)$ is a semishortening derivation in $cf.EQ_{12}$ of 

$\Gamma, s=r \seq F\{v/s\}$.
\ei
\end{proposition}

{\bf Proof} To be more accurate,  ${\cal G}_1$ and ${\cal G}_2$ actually have four arguments, i.e. ${\cal D}$, $F$, $\{v/r\}$ and $s$ and 
their definition requires that $F\{v/r\}$ coincides with the succedent of the  endsequent of ${\cal D}$. However, since it will be clear from the context what $F$ and $\{v/r\}$ are, there is no harm in using the  simplified notations ${\cal G}_1({\cal D}, r,s)$ and ${\cal G}_2({\cal D}, r,s)$.

 By Lemma \ref{singleton}, it suffices to deal with derivations in $cf.EQ_{12}^1$.
If $r\prec s$ then ${\cal G}_1({\cal D},r,s) $ is obtained by applying to ${\cal D}$ an $=_1$-inference
with operating equality  $r=s$ and if $s\not\prec r$ (in particular if $r\prec s$), ${\cal G}_2({\cal D},r,s) $ is obtained by applying to ${\cal D}$ an $=_2$ inference, with operating equality $s=r$.
Hence in defining ${\cal G}_1$ we may assume that $r\not\prec s$, while in defining ${\cal G}_2$ we may assume that $s\prec r$.

${\cal G}_1({\cal D}, r,s)$ and ${\cal G}_2({\cal D},r, s)$   are defined simultaneously  by recursion on
the height $h({\cal D})$ of ${\cal D}$,  for arbitrary $s$.

If $h({\cal D})=0$  we have the following cases.

Case 0.1 ${\cal D}$ reduces to $F\{v/r\}\seq F\{v/r\}$. As ${\cal G}_1({\cal D},r, s)$ we can take

\[
\ba{cl}
F\{v/s\}\seq F\{v/s\}&\vbox to 0pt{\hbox{$=_2^l$}}\\
\cline{1-1}
F\{v/r\},r=s \seq F\{v/s\}&
\ea
\]
which is nonlengthening, since we are assuming that $r\not \prec s$,
and as ${\cal G}_2({\cal D},r, s)$ we can take:
\[
\ba{cl}
F\{v/s\}\seq F\{v/s\}&\vbox to 0pt{\hbox{$=_1^l$}}\\
\cline{1-1}
F\{v/r\}, s=r \seq F\{v/s\}&
\ea
\]
which is shortening, since we are assuming that $s\prec r$.
Thus in both cases we have obtained a semishortening derivation, as required.

Case 0.2 ${\cal D}$ reduces to $~~\seq t_0= t\{v/r\}$  with $t_0\equiv t\{v/r\}$. As ${\cal G}_1({\cal D},r,s)$ 
we can take:

\[
\ba{cl}
\seq t\{v/s\}=t\{v/s\}&\vbox to 0pt{\hbox{$=_2$}} \\
\cline{1-1}
r=s\seq t\{v/r\}=t\{v/s\}
\ea
\]
which is nonlengthening, and as ${\cal G}_2({\cal D},r, s)$ we can take:

\[
\ba{cl}
\seq t\{v/s\}=t\{v/s\}&\vbox to 0pt{\hbox{$=_1$}} \\
\cline{1-1}
s=r\seq t\{v/r\}=t\{v/s\}
\ea
\]
which is shortening.

Case 0.3 ${\cal D}$ reduces to $~~\seq t\{v/r\}= t_0$ with  $t_0\equiv t\{v/r\}$.  The definition of 
${\cal G}_1({\cal D}, r, s)$ and ${\cal G}_2({\cal D}, r, s)$ is essentially the same as in case 0.2.

If $h({\cal D})>0$ and ${\cal D}$  ends with a structural rule and  has the form:

\[
\ba{cl}
{\cal D}_0&\\
\Gamma'\seq F\{v/r\}&\\
\cline{1-1}
\Gamma\seq F\{v/r\}

\ea
\]

 ${\cal G}_1({\cal D},r, s)$ and ${\cal G}_2({\cal D}, r, s)$ are obtained by applying the same structural rule and some exchanges to the endsequent of ${\cal G}_1({\cal D}_0, s)$ and ${\cal G}_2({\cal D}_0, s)$ that, by induction hypothesis, are semishortening  derivations
in $cf.EQ_{12}^1$ of $\Gamma', r=s\seq F\{v/s\}$ and $\Gamma', s=r\seq F\{v/s\}$ respectively.

Otherwise we have the following four cases depending on the ending equality inference of ${\cal D}$.

Case 1. ${\cal D}$ ends with an $=_1$-inference. Then we have the following three subcases:

Case 1.1. ${\cal D}$ has the form:

\[
\ba{cl}
{\cal D}_0&\\
\Gamma' \seq F^\circ\{u/p, v/r\}&\vbox to 0pt{\hbox{$=_1$}}\\
\cline{1-1}
\Gamma', p=q \seq F^\circ\{u/q, v/r\}&
\ea
\]with $F\equiv F^\circ\{u/q\}$ and $v$ does not occurs in $q$.
 Since ${\cal D}$ is semishortening,  $p\prec q$. By induction hypothesis  ${\cal G}_1({\cal D}_0, r,  s)$  is  a semishortening derivation of $\Gamma' , r=s \seq F\{u/p, v/s\}$ and we can let  ${\cal G}_1({\cal D}, r, s)$ be:

\[
\ba{cl}
{\cal G}_1({\cal D}_0,s)&\\
\Gamma' , r=s \seq F^\circ\{u/p, v/s \}&\vbox to 0pt{\hbox{$=_1$}}\\
\cline{1-1}
\Gamma, p=q, r=s \seq F^\circ\{u/q, v/s\}
\ea
\]
The definition of ${\cal G}_2({\cal D}, r,  s)$ is  the same,  except that 
${\cal G}_1({\cal D}_0,r, s)$ and $r=s$ in the endsequent  are replaced by ${\cal G}_2({\cal D}_0, r, s)$ and $s=r$ respectively.

Case 1.2  ${\cal D}$ has the form:

\[
\ba{cl}
{\cal D}_0&\\
\Gamma' \seq F^\circ\{u/p\}&\vbox to 0pt{\hbox{$=_1$}}\\
\cline{1-1}
\Gamma', p=q\{v/r\} \seq F^\circ\{u/q\{ v/r\}\}&
\ea
\]with $F\equiv F^\circ\{u/q\}$ and $v$ occurs in $q$.
By induction hypothesis there is a semishorthening derivation ${\cal G}_1({\cal D}_0, p, q\{v/s\}) $ of 
$\Gamma', p=q\{v/s\} \seq F^\circ\{u/q\{v/s\}\}$ and we can let  ${\cal G}_1({\cal D}, r, s) $ be:

\[
\ba{cl}
{\cal G}_1({\cal D}_0, p, q\{v/s\})&\\
\Gamma', p=q\{v/s\} \seq F^\circ\{u/q\{v/s\}\}&\vbox to 0pt{\hbox{$=_2^l$}}\\
\cline{1-1}
\Gamma', p=q\{v/r\}, r=s \seq F^\circ\{u/q\{v/s\}\}
\ea
\]
which is semishortening,  since its  ending $=_2^l$-inference is nonlengthening, given that in defining ${\cal G}_1$, we are assuming that $r\not\prec s$.

The definition of ${\cal G}_2({\cal D}, r, s)$ is  the same,  except that   $=_2^l$ and $r=s$ in the endsequent
are replaced by  $=_1^l$ and $s=r$ respectively. In fact the ending $=_1^l$-inference of the derivation so obtained is shortening since, in defining ${\cal G}_2$,  we are assuming that $s\prec r$.
 Notice that in this case the definition of ${\cal G}_2({\cal D}, r, s)$ depends on 
${\cal G}_1({\cal D}_0, p, q\{v/s\})$.

Case 1.3 $r$ has the form $r^\circ\{u/q\}$ and  ${\cal D}$  the form:

\[
\ba{cl}
{\cal D}_0&\\
\Gamma' \seq F\{v/r^\circ\{u/p\}\}&\vbox to 0pt{\hbox{$=_1$}}\\
\cline{1-1}
\Gamma', p=q \seq F\{v/r^\circ\{u/q\}\}
\ea
\]
with $p \prec q$.
By induction hypothesis there is a semishortening derivation ${\cal G}_1({\cal D}_0, r^\circ\{u/p\}, s)$ of 
$\Gamma', r^\circ\{u/p\} = s \seq F\{v/s\}$ and we can let ${\cal G}_1({\cal D}, s)$ be:
\[
\ba{cl}
{\cal G}_1({\cal D}_0, r^\circ\{u/p\},  s)&\\
\Gamma', r^\circ\{u/p\} =  s \seq F\{v/s\}&\vbox to 0pt{\hbox{$=_1^l$}}\\
\cline{1-1}
\Gamma', r^\circ\{u/q\}= s, p=q \seq F\{v/s\}
\ea
\]

The definition of ${\cal G}_2({\cal D}, r,  s)$ is  the same,   except that  ${\cal G}_1({\cal D}_0, r^\circ\{u/p\}, s)$, 
$r^\circ\{u/p\}=s $ and $r^\circ\{u/q\}= s$  are replaced by 
${\cal G}_2({\cal D}_0,r^\circ\{u/p\},  s)$, 
$s=r^\circ\{u/p\}$ and $s=r^\circ\{u/q\}$ respectively.

Case 2  ${\cal D}$ ends with an $=_2$ inference.

Case 2.1 ${\cal D}$ has the form:

\[
\ba{cl}
{\cal D}_0&\\
\Gamma' \seq F^\circ\{u/p, v/r\}&\vbox to 0pt{\hbox{$=_2$}}\\
\cline{1-1}
\Gamma', q=p \seq F^\circ\{u/q, v/r\}&
\ea
\]with $F\equiv F^\circ\{u/q\}$ and $v$ does not occur in $q$.
 Since ${\cal D}$ is semishortening,
 $q\not\prec p$.
 By induction hypothesis we have a semishortening derivation ${\cal G}_1({\cal D}_0, r,  s)$ of $\Gamma' , r=s \seq F^\circ\{u/p, v/s\}$ and as ${\cal G}_1({\cal D}, r, s)$ we can take:

\[
\ba{cl}
{\cal G}_1({\cal D}_0, r, s)&\\
\Gamma', r=s \seq F^\circ\{u/p, v/s\}&\vbox to 0pt{\hbox{$=_2$}}\\
\cline{1-1}
\Gamma, q=p , r=s \seq F^\circ\{u/q, v/s\}
\ea
\]

The definition of ${\cal G}_2({\cal D}, r,  s)$ is  the same,   except that    ${\cal G}_1({\cal D}_0, r, s)$ and $r=s$ in the endsequent are  replaced by ${\cal G}_2({\cal D}_0, r, s)$ and $s=r$ respectively.

Case 2.2  ${\cal D}$ has the form:

\[
\ba{cl}
{\cal D}_0&\\
\Gamma' \seq F^\circ\{u/p\}&\vbox to 0pt{\hbox{$=_2$}}\\
\cline{1-1}
\Gamma', q\{v/r\}=p \seq F^\circ\{u/q\{ v/r\}\}&
\ea
\]with $F\equiv F^\circ\{u/q\}$ and $v$ occurs in $q$.
By induction hypothesis there is a semishortening derivation ${\cal G}_2({\cal D}_0, p,  q\{v/s\}) $ of 
$\Gamma', q\{v/s\}= p \seq F^\circ\{u/q\{v/s\}\}$ and we can let  ${\cal G}_1({\cal D}, r,  s) $ be:

\[
\ba{cl}
{\cal G}_2({\cal D}_0, p, q\{v/s\})&\\
\Gamma', q\{v/s\}=p \seq F^\circ\{u/q\{v/s\}\}&\vbox to 0pt{\hbox{$=_2^l$}}\\
\cline{1-1}
\Gamma', q\{v/r\} = p, r=s \seq F^\circ\{u/q\{v/s\}\}
\ea
\]
which  is semishortening,  since its  ending $=_2^l$-inference is nonlengthening, given that in defining ${\cal G}_1$, we are assuming that $r\not\prec s$.  Notice that in this case the definition of ${\cal G}_1({\cal D}, r,  s)$ depends on 
${\cal G}_2({\cal D}_0, p, q\{v/s\})$.

The definition of ${\cal G}_2({\cal D},  r, s)$ is the same,  except that $=_2^l$ and $r=s$ in the endsequent are replaced by   $=_1^l$ and $s=r$ respectively .

Case 2.3 $r$ has the form $r^\circ\{u/q\}$ and  ${\cal D}$  the form:

\[
\ba{cl}
{\cal D}_0&\\
\Gamma' \seq F\{v/r^\circ\{u/p\}\}&\vbox to 0pt{\hbox{$=_2$}}\\
\cline{1-1}
\Gamma', q=p \seq F\{v/r^\circ\{u/q\}\}
\ea
\]
with $q \not \prec p$.
By induction hypothesis,   ${\cal G}_1({\cal D}_0, r^\circ\{u/p\},  s)$ is  a semishorthening derivation    of 
$\Gamma', r^\circ\{u/p\} = s \seq F\{v/s\}$ and we can let ${\cal G}_1({\cal D}, r, s)$ be:
\[
\ba{cl}
{\cal G}_1({\cal D}_0, r^\circ\{u/p\}, s)&\\
\Gamma', r^\circ\{u/p\} =  s \seq F\{v/s\}&\vbox to 0pt{\hbox{$=_2^l$}}\\
\cline{1-1}
\Gamma',  q=p, r^\circ\{u/q\}= s \seq F\{v/s\}
\ea
\]
The definition of ${\cal G}_2({\cal D}, r, s)$ is  the same, except that  ${\cal G}_1({\cal D}_0, r^\circ\{u/p\}, s)$  and $r^\circ\{u/q\}=s$ are replaced by  
${\cal G}_2({\cal D}_0, r^\circ\{u/p\}, s)$ and $s=r^\circ\{u/q\}$  respectively.

Case $1^l$  ${\cal D}$ ends with an $=_1^l$-inference, i.e.  it has the form:

\[
\ba{cl}
{\cal D}_0&\\
\Gamma', G\{u/p\}\seq F\{v/r\}&\vbox to 0pt{\hbox{$=_1^l$}}\\
\cline{1-1}
\Gamma', G\{u/q\}, p= q\seq F\{v/r\}
\ea
\]
with $p\prec q$. By induction hypothesis  ${\cal G}_1({\cal D}_0, r,  s)$ is a semishortening derivation
of
$\Gamma', G\{u/p\}, r=s \seq F\{v/s\}$ and we can let ${\cal G}_1({\cal D}, r,  s)$ be:

\[
\ba{cl}
{\cal G}_1({\cal D}_0, r, s)&\\
\Gamma', G\{u/p\}, r=s \seq F\{v/s\}&\vbox to 0pt{\hbox{$=_1^l$}}\\
\cline{1-1}
\Gamma', G\{u/q\}, p=q, r=s \seq F\{v/s\}
\ea
\]
The definition of ${\cal G}_2({\cal D}, r, s)$ is the same,  except that  ${\cal G}_1({\cal D}_0, r, s)$ and $r=s$ in the endsequent are replaced  by 
${\cal G}_2({\cal D}_0, r, s)$ and $s=r$ respectively.

Case $2^l$  ${\cal D}$ ends with a $=_2^l$-inference, namely $p=q$ is replaced by $q=p$ in the endsequent of ${\cal D}$ as represented in  Case $1^l$. Then ${\cal G}_1({\cal D}, r, s)$ and ${\cal G}_2({\cal D}, r, s)$ are defined as 
in Case $1^l$,  except that $p=q$ is replaced by $q=p$.
$\Box$

\begin{theorem} \label{semishortening}
Any derivation in $EQ_{12}$ can be tranformed into a cut-free semishortening derivation in $EQ_{12}$ of its endsequent.

\end{theorem}

{\bf Proof} Every derivation in $EQ_{12}$ can be effectively transformed into a derivation in $EQ$, henceforth, by Theorem \ref{CutEQ}, into a cut free derivation in $EQ$ of its endsequent. The conclusion follows by the admissibility of the equality rules $=_1$ and $=_2$ in $cfEQ_{12}$ restricted to semishortening derivations, established in the previous Proposition \ref{gamma}. $\Box$

\begin{corollary}
Any derivation in $LJ_{12}$ or $LK_{12}$ can be tranformed into a cut-free semishortening derivation in the same calculus  of its endsequent. 
\end{corollary}

{\bf Remark} Since the semishortening derivations are nonlengthening,  an immediate consequence of  Theorem  \ref{semishortening} is that every derivation  in $EQ_{12}$ can be transformed into a cut-free nonlengthening derivation of its endsequent.
That can also be established  by observing that if   semishortening is replaced by nonlengthening, then 
Proposition \ref{gamma} still hold with 
essentially the same proof. 

\ 

{\bf Remark}  Since  $=_1^l$ and $=_2^l$ are derivable in $EQ$, from the admissibility of the cut rule in $cf.EQ$, it follows that  $=_1^l$ and $=_2^l$ are admissible in $cf.EQ$. Hence from
 Proposition \ref{gamma} 
it follows that also  $=_1^l$ and $=_2^l$ are admissibile in $cf.EQ_{12}$ restricted to semishortening derivations. As it results from \cite{L68}, in the case of nonlengthening derivations,  a   direct inductive proof of this admissibility result is possible,
 but it requires the additional assumption that $\prec$ be a strict partial order congruent with respect to substitution. It is the admissibility of $=_1^l$ and $=_2^l$ in $cf.EQ$, unnoticed in \cite{L68}, that allows for the weakening of such assumption to the requirement that $\prec$ be simply  antisymmetric.

\subsection{Related work}

Beside the references already given in the introduction, we add that the restriction $=_0\seq$ of the rule $=\seq$ to atomic $F$ had been considered in \cite{Na66} in conjunction with the following {\em left reflexivity elimination rule}:
\[
\ba{c}
\Gamma, t=t\seq \Theta\\
\cline{1-1}
\Gamma\seq \Theta
\ea
\]
in the framework of $LJ$ and $LK$.  We point out that for the resulting systems, cut  elimination is a trivial matter as any formula $H$  can be seen (in many ways) as $H^\circ\{v/t\}$, so that the cut rule can be derived from $=\seq$ and the left reflexivity elimination rule, as follows:
\[
\ba{c}
\Gamma \seq \Delta, H~~~~~\Lambda, H \seq \Theta\\
\cline{1-1}
\Gamma, \Lambda, t=t \seq \Delta, \Theta\\
\cline{1-1}
\Gamma, \Lambda\seq \Delta, \Theta\\
\ea
\]
 Despite the fact that  the left reflexivity elimination rule eliminates equalities, \cite{Na66} uses such systems  to prove the conservativity of first order logic with equality over first order logic without equality.
  A calculus similar to $LK_{N}^{=}$, but without $\vel$ and $\exists$ and with 
$CNG$ restricted from the start to atomic formulae is considered in \cite{Pl71},  which establishes the eliminabilty of the cut rule without going through the reduction of derivation to separated form.
The idea of using the admissibility of the rule $CNG$ in $EQ_N$ to prove Theorem \ref{cutL} first appeared in 
\cite{PP13}. However the proof of  admissibility and the way of deriving the cut elimination theorem for the $LJ^=$ and $LK^=$ systems  given in this paper are a substantial improvement of those to be found in \cite{PP13}.


\begin{thebibliography}{10}


\bibitem{C77} 
H.B. Curry,
\newblock Foundations of Mathematical Logic
\newblock Dover (1977)



\bibitem{G35}
G. Gentzen, 
\newblock{ Untersuchungen uber der logische Schliessen}
\newblock{\em Matematische Zeitschrift}
Vol.39,  pp. 176-210, 405-431 (1935)


\bibitem{Ga86}
J.  Gallier, 
\newblock  Logic for Computer Science, Foundations of Automatic Theorem Proving,
Harper \& Row,  NewYork (1986).


\bibitem{K63}
S.  Kanger,
\newblock A Simplified Proof  Method for Elementary Logic.
\newblock  In: P. Braffort, D. Hirshberg (eds)
\newblock Computer Programming and Formal Systems, pp. 87-94.
 North-Holland, Amsterdam (1963)


\bibitem{L68}
A.V.  Lifschitz,
\newblock  Specialization of the form of deduction in the predicate calculus  with equality and function symbols (in Russian). In Trudy MIAN, vol. 98, 5-25 (1968). 
\newblock English translation in: V.P. Orevkov  (ed) {\em The Calculi of Symbolic Logic. I}, Proceedings of the Steklov Institute of Mathematics 98 (1971)




\bibitem{Na66}
T.  Nagashima, 
\newblock An Extension of the Craig-Schutte Interpolation Theorem
\newblock  Annals of the Japan Association for the Philosophy of Science 3, 12-18, (1966)




\bibitem{NvP98} J. von Plato, S.  Negri,
\newblock  Cut Elimination in the Presence of Axioms. 
\newblock The Bulletin of Symbolic  Logic 4 (4) , 418--435 (1998)


\bibitem{O69}
V. P. Orevkov,
\newblock  On Nonlengthening Applications of Equality Rules (in Russian)
\newblock Zapiski Nauchnyh Seminarov LOMI, 16:152-156, 1969
\newblock English translation in: A.O. Slisenko  (ed) {\em Studies in Constructive Logic}, Seminars in Mathematics:  Steklov Math. Inst. 16, Consultants Bureau, NY-London 77-79 (1971)


\bibitem{vP12} J. von Plato, 
\newblock  Gentzen's Proof Systems: Byproducts in a Work of Genius
\newblock The Bulletin of Symbolic  Logic 18 (3) , 317- 367 (2012)


\bibitem{vP14} J. von Plato, 
\newblock  From Axionatic Logic to Natural Deduction
\newblock Studia Logica 102 (6) , 1167-1184 (2014)


\bibitem{NvP01}
J. von Plato, S.  Negri,
\newblock  Structural Proof Theory.
 \newblock Cambridge University Press  (2001)



\bibitem{Pr65}
D. Prawitz,
\newblock Natural Deduction. A Proof-Theoretical  Study
\newblock Almquist and Wiksell  (1965)


\bibitem{PP13}
F. Parlamento, F.  Previale, 
\newblock Cut elimination for Gentzen's Sequent Calculus with Equality and Logic of Partial Terms.
\newblock Lecture Notes in Computer Science 7750, 161-172 (2013)



\bibitem{Pl71}
R.A.  Pliuskevicius, 
\newblock  A sequential Variant of Contructive Logic Calculi for Normal Formulas not Containing Structural Ruels.
\newblock In  {\en The Calculi of Symbolivc Logic 1} V. Orenkov  ed. Proceeding of the Steklov Institute  of Mathematics 98. 175-229 (1971)


\bibitem{S69}
  M.E. Szabo ed,
\newblock  The Collected Papers of Gerhard Gentzen,
\newblock North Holland,
(1969)




\bibitem{T87} 
G. Takeuti,
\newblock  Proof Theory. Studies in Logic and the Foundations of Mathematics vol 81, 2nd edition
\newblock  North Holland, Amsterdam (1987)

\bibitem{TS96}
A.S.  Troelstra, H.  Schwichtemberg, 
\newblock   Basic Proof Theory. Cambridge Tracts in Theoretical Computer Science, Vol. 43
Cambridge University Press, Cambridge(1996).


\bibitem{TS00} 

A.S. Troelstra, H.  Schwichtemberg,
\newblock    Basic Proof Theory. 2nd edition
\newblock Cambridge University Press, Cambridge(2000).






\end{thebibliography}
 \end{document}